\documentclass[12pt]{article}
\usepackage[top=1in,bottom=1in, left=1in, right=1in]{geometry}
\usepackage[utf8]{inputenc}
\usepackage{bbm}
\usepackage{mathtools}
\usepackage{enumitem}
\usepackage{diagbox}
\usepackage[english]{babel}
\usepackage{amssymb}
\usepackage[sort,compress]{natbib}
\usepackage{amsmath}
\usepackage{xcolor}
\usepackage{url}
\usepackage{relsize}
\usepackage{breqn}
\usepackage{comment}
\usepackage{amsmath,amsthm,amssymb,amsopn,amsfonts,pdfpages,dsfont}
\usepackage{graphics}
\usepackage{graphicx}
\usepackage{bbm,subfig}
\usepackage{soul}
\usepackage{extarrows}
\usepackage{caption,url}
\usepackage[colorlinks=true,linkcolor=cyan,citecolor=blue]{hyperref}
\usepackage{booktabs}
\usepackage{authblk}
\usepackage{bbold}
\usepackage{setspace}
\usepackage{blkarray}
\usepackage{multirow}
\usepackage{float}
\usepackage{fancyvrb}
\usepackage{etoolbox}
\usepackage{cancel}
\patchcmd{\thmhead}{(#3)}{#3}{}{}

\newtheorem{Def}{Definition}
\newtheorem*{Def*}{Definition}
\newtheorem{lem}{Lemma}
\newtheorem{rem}{Remark}
\newtheorem{thm}{Theorem}

\newcommand{\E}{\mathbb{E}}
\newcommand{\Prb}{\mathbb{P}}
\newcommand{\PP}{\mathbb{P}}

\newcommand{\ppt}{\mathfrak{p}_t}

\newcommand{\aij}{A_{i,j}}
\newcommand{\snp}{S_{n,P}}
\newcommand{\bij}{B_{i,j}}
\newcommand{\bsij}{B_{\sigma_P(i),\sigma_P(j)}}
\newcommand{\tr}[1]{\operatorname{Tr}\!\left(#1\right)}

\newcommand{\ER}{\operatorname{ER}}
\newcommand{\SBM}{\operatorname{SBM}}
\newcommand{\tcov}{t_{\mathrm{cover}}}

\newcommand{\Bern}{\operatorname{Bern}}
\newcommand{\Bin}{\operatorname{Bin}}
\newcommand{\Geo}{\operatorname{Geo}}
\newcommand{\Corr}{\operatorname{Corr}}

\newcommand{\trel}{t_{\text{rel}}}

\DeclareMathOperator*{\argmax}{arg\,max}
\DeclareMathOperator*{\argmin}{arg\,min}

\newcommand{\trz}{\textcolor{black}{0}}
\title{Matching and mixing: Matchability of graphs under Markovian error}
\author[$1$]{Zhirui Li}
\author[$2$]{Keith D. Levin}
\author[$3$]{Zhiang Zhao}
\author[$4$]{Vince~Lyzinski}

\affil[$1$]{\small Department of Statistics and Data Science, University of Pennsylvania}
\affil[$2$]{\small Department of Statistics, University of Wisconsin--Madison}
\affil[$3$]{\small Department of Electrical and Computer Engineering, University of Kentucky}
\affil[$4$]{\small Department of Mathematics, University of Maryland, College Park}

\begin{document}
\maketitle
\begin{abstract}
    We consider the problem of graph matching for a sequence of graphs generated under a time-dependent Markov chain noise model. 
    Our edgelighter error model, a variant of the classical lamplighter random walk, iteratively corrupts the graph $G_0$ with edge-dependent noise, creating a sequence of noisy graph copies $(G_t)$.
    Much of the graph matching literature is focused on anonymization thresholds in edge-independent noise settings, and we establish novel anonymization thresholds in this edge-dependent noise setting when matching $G_0$ and $G_t$. 
    Moreover, we also compare this anonymization threshold with the mixing properties of the Markov chain noise model.
    We show that when $G_0$ is drawn from an Erd\H{o}s–R\'{e}nyi model, the graph matching 
    anonymization threshold and the mixing time of the edgelighter walk are both of order $\Theta(n^2\log n)$. 
    We further demonstrate that for more structured model for $G_0$ (e.g., the Stochastic Block Model), graph matching anonymization can occur in $O(n^\alpha\log n)$ time for some $\alpha<2$, indicating that anonymization can occur before the Markov chain noise model globally mixes.
    Through extensive simulations, we verify our theoretical bounds in the settings of Erd\H{o}s–R\'{e}nyi random graphs and stochastic block model random graphs, and explore our findings on real-world datasets derived from a Facebook friendship network and a European research institution email communication network. 
\end{abstract}

\section{Introduction}

Networks are used across the sciences to describe complex real-world relationships.
Examples include studies of social networks, where nodes represent users and edges denote friendships or follower-following relations \cite{mcauley2012learning}. 
Co-authorship networks \cite{fonseca2016co} and email exchange networks \cite{leskovec2007graph} have proven to be helpful in identifying leading research scientists, prominent organizations, and potential academic fraud.
In neuroscience, networks derived from functional MRI scans \cite{wang2010graph,Sporns2012} have seen applications in analyzing connectivity patterns \cite{farahani2019application,mishchencko2011bayesian}, and studying brain-related diseases \cite{wang2013disrupted, vecchio2017connectome}.
Large knowledge graphs have become increasingly useful for representing the conceptual relations between entities \cite{zou2020survey}.

Statistical network inference on collections of graphs often assumes that the observed graphs are noisy copies generated from some shared latent structures \cite{wang2019joint, levin2017central,priebe2015statistical,arroyo2021maximum,LLL2022}. 
From these observations, numerous tasks such as classification, hypothesis testing, and more can be pursued \cite{vogelstein2011shuffled, relion2019network, MT2, du2023hypothesis,ginestet2017hypothesis,cape2019two,tang2013universally,de2010advantages}.
A crucial assumption of most of these methods is that the vertices are {\em a priori} aligned across networks before inference is pursued.
In situations where this is not the case, inference can degrade when vertex labels are shuffled or noisily observed \cite{saxena2025lost}.
In these cases, graph matching methods are often deployed to attempt to recover the true node labeling across graphs and enable the application of subsequent inference methods \cite{Lazaro2025, wu2023testing, alnaimy2022expanded, frigo2021network, arroyo2021maximum}. 

Stated simply, the goal of graph matching is to find the best possible alignment between the vertices of two (or more) graphs in order to minimize a measure of edgewise structural discrepancy.
For detailed discussions on the graph matching problem, including its applications, approximation methods, extensions, and variants, refer to the surveys \cite{riesen2010exact,yan2016short,conte2004thirty}.
Recently, there have been numerous papers in the literature dedicated to algorithmic graph matching \cite{FAQ,gm1,gm2,gm3,gm4,gm5} and the information theoretic limits of matching theory \cite{Gross2013PGM,lyzinski2016information,cullina2016improved,cullina2017exact,wu2021settling}.
Nearly all of the aforementioned theoretical graph matching developments rely on the assumption of independent edge-level noise.
However, in real applications such as protein-protein interaction networks \cite{stelzl2005human}, epidemiological networks \cite{Newman2002Disease}, error propagation in neural networks \cite{li2017understanding}, and temporal networks \cite{athreya2024euclideanmirrorsdynamicsnetwork}, it is more realistic to consider time- or node-dependent noise.
To address this issue, our focus in the present work concerns the matchability of graphs with spatially- or temporally-dependent edge-noise processes. 

\subsection{Notation}
\label{sec:not}

For any $n\in \mathbb{Z}^+$, we write $[n]= \{1,2,\dots, n\}$.
For $x \in \mathbb{R}$, we write $\lfloor x\rfloor$ for the floor function of $x$, i.e., the element $n\in \mathbb{Z}$ such that $n\leq x$ and $n+1>x$ and denote the ceiling function analogously by $\lceil\cdot\rceil$.
For a set $S$, we write $\binom{S}{j}$ to denote the collection of all size-$j$ subsets of $S$; 
we also write $\{0,1\}^{\binom{S}{j}}$ to denote the set of length $|\binom{S}{j}|$ binary sequences indexed by $\binom{S}{j}$.
We use $\|\cdot\|_F$ to denote the matrix Frobenius norm.
We use the standard notations for asymptotics.
In particular, for two functions $f, g: \mathbb{Z}^+\rightarrow \mathbb{R}^+$, we write $f(n)=o(g(n))$ or $f\ll g$, if $\lim_{n\rightarrow\infty} f(n)/g(n)=0$;
we write $f(n)=\omega(g(n))$ or $f\gg g$ if $g(n)= o(f(n))$; 
we write $f(n)=O(g(n))$ or $f\lesssim g$ if $\exists C>0$ and $n_0$ such that $\forall n\geq n_0$, $f(n)\leq Cg(n)$;
we write $f(n)=\Omega(g(n))$ or $f\gtrsim g$ if $g(n) = O(f(n))$;
we write $f(n)=\Theta(g(n))$ if $f(n)=O(g(n))$ and $g(n)=O(f(n))$.

\subsection{The Graph Matching Problem}
\label{sec:GMP}

The formulation of the graph matching problem we consider herein is defined as follows.
\begin{Def} \label{def:GMP}
  Consider two graphs $G_1=(V_1,E_1)$ and $G_2 = (V_2,E_2)$ with corresponding adjacency matrices $A$ and $B$, respectively.
  When $|V_1|=|V_2|=n$, the Graph Matching Problem (GMP) seeks to find
  \begin{equation} 
   \label{eq:GMP:objective_TraceForm}
  \argmin_{P\in\Pi_n} \| A-PBP^T \|_F=
  \argmax_{P\in\Pi_n} \tr{APBP^T},
  \end{equation}
  where $\Pi_n$ denotes the set of all $n$-by-$n$ permutation matrices.
\end{Def} 

\noindent The goal of the GMP is to find the best possible alignment of the vertices (as encoded by a permutation of the vertex labels) between the two graphs, so as to minimize a measure (here the $\ell_2$ norm of binary edge-difference vectors) of the edgewise structural differences. 
In the general setting where the two graphs are allowed to be weighted, directed and to exhibit self-loops, the above formulation is equivalent to the NP-hard quadratic assignment problem \cite{gm_relax, VogConLyzPodKraHarFisVogPri2015}. 
That said, there are efficient graph matching algorithms for special classes of graphs (e.g., planar graphs) \cite{schmidt2009efficient, anari2020planar}. 
When discussing the GMP in what follows, we will refer to a graph and its adjacency matrix interchangeably, due to the one-to-one correspondence between weighted graphs and square matrices.

The GMP has a long and rich history. Initial research focused on finding graph isomorphisms, solutions for which the objective function in Equation~\eqref{eq:GMP:objective_TraceForm} equals zero \cite{babai2015graph,read1977graph, takapoui2016linear}. 
The stringent nature of exact matching, however, limits the applicability of the graph isomorphism problem in real-world scenarios where noise or distortion of graph structures might be present. In response, more general formulations of the graph matching problem aim not to drive the objective to zero, but rather to minimize it, allowing for the possibility that two graphs, while not precisely isomorphic, may still exhibit meaningful structural similarity; see \cite{conte2004thirty} for more extensive background.

In the graph isomorphism setting, there is a notion of a true alignment (i.e., the isomorphism) that is not necessarily present in the general GMP.
To remedy this, we consider the pair of networks $G_1$ and $G_2$ as having a true, latent (i.e., a priori unknown) node-level alignment across their vertex sets.
This alignment could be data driven: In a social network, it could represent an account belonging to the same user across networks;
in connectomics, it could represent the same neuronal region from different scans of a single subject; 
in protein-protein interaction networks, this could represent the same protein in two different interaction networks.
This alignment could also be model driven, as we will see in Section \ref{sec:ER}, where our random graph models have correlation across edges in $G_1$ and $G_2$ which identifies vertices across graphs.
Herein, to specify a ground truth alignment, we use the terms \emph{node-aligned} or \emph{vertex-aligned} to refer to a pair of graphs $G_1$ and $G_2$ (with resp., adjacency matrices $A$ and $B$) and a latent alignment matrix $P$ (wlog, we often set $P$ to be the identity alignment $I_n$) such that for all $i \in [n]$, node index $i$ in $A$ and node index $\sigma_P(i)$ in $B$ (where $\sigma_P$ is the permutation function associated with the permutation matrix $P$) represent the same entity in each graph.
Irrespective of algorithm, we then define
the term {\em matchability} of graphs as the information-theoretic possibility of the GMP to recover a ground truth network alignment; i.e., if the latent alignment is $P^*$, is it the case that $\{P^*\}=\argmax_{P\in\Pi_n} \tr{APBP^T}$. 

In suitable random graph models, the literature has identified scenarios where matchability is theoretically achievable \cite{wu2021settling} and practically achieved \cite{mossel2020seeded,fan2020spectral}.
As mentioned above, often the node-alignment across graphs is determined by introducing correlation between the edges across the networks, 
and recent work has identified a sharp threshold on this edge correlation for graph matchability to be feasible \cite{wu2021settling,lyzinski2016information,cullina2016improved,cullina2017exact}.
In these models, edges within each network $G_1$ and $G_2$ are often assumed to be independent, the edges across $G_1$ and $G_2$ are \emph{$\rho$-correlated} as follows.
For each pair of vertices $\{u,v\}\in\binom{V}{2}$, we have 
$$\text{corr}(\mathds{1}\{\{u,v\}\in E(G_1)\},\mathds{1}\{\{u,v\}\in E(G_2)\})=\rho.$$ 
Moreover, it is often assumed that the collection of variables $$\bigg\{\mathds{1}\{\{u,v\}\in E(G_1)\},\mathds{1}\{\{w,z\}\in E(G_2)\}\bigg\}_{\{u,v\},\{w,z\}\in\binom{V}{2}}$$ are independent outside of this correlation.
Stated informally, in moderately dense graphs a correlation of order at least $\sqrt{\log n/n} $ is required for matchability with high probability. 
In the setting of these sharp threshold results, distinct pairs of edges (e.g., $\{u,v\}$ and $\{w,z\}$) are independent, both within a single graph and across graphs.
In cases where a dependent noise structure exists, as explored herein, only a few results have been established to formulate analogous threshold results; see, for example, \cite{arroyo2021maximum}.

\subsection{Random Graph Models} \label{sec:models}
To study the behavior of graphs and perform statistical inference, researchers often restrict themselves to low-rank graph models (i.e., where $\mathbb{E}(A)$ is low-rank). 
Among these models, the Erd\H{o}s–R\'{e}nyi (ER) graph model \cite{Erdos} is the simplest, as it assumes that all edges are present or absent independently with the same probability.
\begin{Def} \label{def:er}
   Let $n$ be a positive integer and $p\in[0,1]$.
   A random graph $G$ on $n$ vertices is distributed according to an Erd\H{o}s–R\'{e}nyi Model (ER) with parameters $(n,p)$, denoted as $G \sim \ER(n,p)$, if for any pair of vertices $i,j\in V(G)$ with $i\neq j$, we have that an edge exists in $G$ between $i$ and $j$ with probability $p$, independent of all other edges.
\end{Def}
\noindent Though simple, the Erd\H os-R\'enyi model has proven to be a rich model for mathematical analysis, with various thresholds for graph properties established under this model \cite{FK16,bollobas1998random,Erdos}.
See also \cite{bollobas07} for the ``inhomogeneous'' Erd\H os-R\'enyi graph, in which edges are present or absent independently but with possibly different probabilities.

While our first result on the interplay between the mixing of our Markov (dependent) noise process and graph matchability is set in the Erd\H{o}s-R\'{e}nyi model, more nuanced results require graph models with more nuanced structure. 
To this end, we consider an  extension of the Erd\H os-R\'enyi model that posits community structured edge connectivity: the stochastic block model (SBM) \cite{sbm}. 
\begin{Def} \label{def:sbm}
Let $n$ and $K$ be positive integers with $n \ge K$ and let $\tau : [n] \rightarrow [K]$.
Given a symmetric matrix $\Lambda \in [0,1]^{K \times K}$, a graph on $n$ vertices is distributed according to a Stochastic Block Model (SBM) with parameters $K, \Lambda, \tau$, denoted $G \sim \SBM(K,\Lambda, \tau)$, if for any pair of nodes $u,v$, we have 
    $$\mathbb{P}(\text{there is an edge between }u\text{ and } v)=\Lambda_{\tau(u),\tau(v)}.$$
    Moreover, the collection of events 
    $$\bigg\{ \{\text{there is an edge between }u\text{ and } v\}\bigg\}_{\{u,v\}\in\binom{V}{2}}$$
    are independent conditional on $\tau$.
The vector $\tau$ encodes a partition of the $n$ vertices into $K$ communities (i.e., ``blocks''), which we denote by $C_k = \{ i : \tau(i)=k \}$ for each $k \in [K]$, and we define the vector $\vec{n} = (n_1,n_2,\dots,n_K)$ of community sizes $n_k=|C_k|$.
\end{Def}
\noindent In essence, the SBM model generalizes the ER random graph to allow for edges to differ in their probabilities according to the community memberships of the vertices, rather than having all edges present or absent with the same probability.
The SBM has been very well studied in the literature, as it has proven to be a rich model for theoretically and algorithmically exploring the properties of community detection in networks; see, for example, \cite{lyzinski2016community, abbe2018community, lee2019review,bickel2009nonparametric,rohe2011spectral,STFP,karrer11stoch}.

\subsection{Markov Error Model} \label{sec:lamp}

Given a graph $G$, a random walk on the vertices of $G$ can be defined by having a random walker traverse, at each step, from one vertex to a (possibly neighboring) vertex chosen uniformly at random \cite{lovasz1993random,levin2017markov}.
In the present work, we will construct all of our random walks to be Markov chains with appropriately defined state spaces and transition matrices. 
A key property of Markov chains is that ergodic (i.e., positive recurrent, aperiodic and irreducible) chains converge to their stationary distribution \cite{levin2017markov}.
If the random walker is corrupting the edge structure of a given graph, and the stationary distribution of the chain (defined as a distribution over the space of all networks of appropriate order) is sufficiently flat/uniform, then it stands to reason that as the chain approaches stationarity, the identity of individual vertices is eventually anonymized.  
This idea/interplay forms the basis for the theoretical developments below.

To wit, we construct a lamplighter-like walk in which the lamps are on edges rather than vertices of the graph.
In this walk, the random walker turns the lamp on or off with a certain probability each time it traverses an edge. 
By considering the position of the lighter after each step along with the status of the lamps on all edges, we obtain a Markov Chain on the space $V(G)\times \{1,\trz\}^{\binom{V(G)}{2}}$. 
That is, the state is captured by the current location of the random walker and the status of the lamps at each of the possible edges of $G$ is encoded by an element of $\{1,\trz\}^{\binom{V(G)}{2}}$. 
See \cite{levin2017markov} for a detailed discussion of the similarly constructed lamplighter walks where lamps are on vertices.
The model we consider for introducing error onto the edges of a graph is inspired by these lamplighter walks and is defined as follows.

Let $G=(V,E)$ be a graph and consider the function $h_G:\binom{V}{2}\to\{\trz,1\}$, defined according to 
\begin{equation*}
h_G\left( \{u,v\} \right)
=   \begin{cases}
    1 &\mbox{ if } \{u,v\}\in E \\
    \trz &\mbox{ otherwise. }
    \end{cases}
\end{equation*}
We will say that an edge $e=\{i,j\} \in E$  is ``on'' (i.e., $h(e)=1$) and elements of $e\in\binom{V}{2}\setminus E$ (i.e., for which $h(e)=\trz$) are ``off''.
Note that by switching the ``lamps'' on or off, the edgelighter changes $E$, and thus the graph $G$.
We use $G_0$ to denote the initial graph and $G_t$ to denote the graph after $t$ steps with corresponding adjacency matrices $A_0$ and $A_t$. 

\begin{Def}[(Standard Edgelighter Walk)] \label{def:stdwalk} 
    Let $G_t=(V,E_t)$ be the state of the graph at time $t$ with $G_0=G=(V,E_0=E)$.
    The edgelighter is modeled via a time-homogeneous Markov chain $Z_t$ with state space $V\times \{\trz,1\}^{\binom{V}{2}}$. 
    If at time $t$ the edgelighter is at vertex $u$, then the edgelighter selects a vertex $v\in V$
    uniformly at random and moves to that vertex.
    If the edgelighter selected $v=u$ at time $t+1$, then $G_{t+1}=G_t$.
    If the edgelighter selects a vertex $v\neq u$ at time $t+1$, we consider two cases.
    \begin{enumerate}
    \item[i.] If $h_{G_{t}}(\{u,v\})=1$, then the graph evolves via 
    $G_{t+1}=(V,E_{t+1})$ with 
    \begin{align}
    \label{eq:lamp_on_to_off}
    h_{G_{t+1}}(e)
    =\begin{cases}
    h_{G_{t}}(e) & \mbox{ if }e \neq \{u,v\};\\
        X_t(1-h_{G_{t}}(e))+(1-X_t)h_{G_{t}}(e)& \mbox{ if }e= \{u,v\}
    \end{cases}
    \end{align}
    where $X_t\sim\Bern(q_1)$ is independent of the walker and the graph history $(G_s)_{s\leq t}$.
    \item[ii.] If $h_{G_{t}}(\{u,v\})=\trz$, then the graph evolves via 
    $G_{t+1}=(V,E_{t+1})$ where 
    \begin{align}
    \label{eq:lamp_off_to_on}
    h_{G_{t+1}}(e)=\begin{cases}
    h_{G_{t}}(e)& \text{ if }e\neq \{u,v\};\\
        Y_t(1+h_{G_{t}}(e))+(1-Y_t)h_{G_{t}}(e)& \text{ if }e= \{u,v\}
    \end{cases}
    \end{align}
    where $Y_t\sim \Bern(q_2)$ is  independent of the walker and the graph history $(G_s)_{s\leq t}$. 
    \end{enumerate}
 Stated simply, if the edgelighter traverses an ``on'' edge (respectively, an ``off'' edge), then independently with probability $q_1$ (resp., $q_2$) the edgelighter  switches the edge from on to off (resp., from off to on), and with probability $1-q_1$ (resp., $1-q_2$) leaves the edge on (resp., off).
 \end{Def}
 We note here that the chain defined above is aperiodic, irreducible, and reversible with respect to the measure $\pi^\circ(\,(u,c)\,)\propto \frac{1}{n}q_2^{\#\text{ of on edges in }c}q_1^{\#\text{ of off edges in }c}$. 
The first two assertions are immediate. Reversibility follows from the following argument.
Let $(u,c)$ and $(v,c')$ be two states in $V\times \{1,\trz\}^{\binom{V}{2}}$ such that $c$ and $c'$ are equal on all possible edges except, perhaps, at $\{u,v\}$ (where $u\neq v$).
Note that if $c=c'$, then 
\begin{equation*}
\mathbb{P}_{(u,c)}(Z_1=(v,c'))=\mathbb{P}_{(v,c')}(Z_1=(u,c))=\begin{cases}
\frac{1}{n}(1-q_1)\text{ if }c(\{u,v\})=1;\\
\frac{1}{n}(1-q_2)\text{ if }c(\{u,v\})=\trz.
\end{cases}
\end{equation*}
If $c(\{u,v\})\neq c'(\{u,v\})$, then we have
\begin{itemize}
\item[] {\bf Case 1:} $c(\{u,v\})=1$, in which case
        \begin{equation*}
        \mathbb{P}_{(u,c)}[Z_1=(v,c')]=q_1/n
        \text{ and }
        \mathbb{P}_{(v,c')}[Z_1=(u,c)]=q_2/n.
        \end{equation*}
\item[] {\bf Case 2:} $c(\{u,v\})=\trz$, in which case
        \begin{equation*}
        \mathbb{P}_{(u,c)}[Z_1=(v,c')]=q_2/n
        \text{ and }
        \mathbb{P}_{(v,c')}[Z_1=(u,c)]=q_1/n.
        \end{equation*}
\end{itemize}
Let $m$ be the number of $1$'s in $c$ (i.e., number of edges in $c$) and $r$ the number of $\trz$'s (i.e., non-edges).
In Case 1,
\begin{align*}
\pi^\circ(\,(u,c)\,) \mathbb{P}_{(u,c)}(Z_1=(v,c'))&\propto\frac{1}{n}q_2^{m}q_1^{r}\frac{1}{n}q_1,\\
\pi^\circ(\,(v,c')\,) \mathbb{P}_{(v,c')}(Z_1=(u,c))&\propto\frac{1}{n}q_2^{m-1}q_1^{r+1}\frac{1}{n}q_2,
\end{align*}
and the two are equal.
Case 2 follows analogously and details are omitted.
Ergodicity of $(Z_t)$ then follows immediately.

\begin{rem}
\label{rem:global}
\emph{
If we let $L_t$ denote the position of the edgelighter random walk at time $t$, then the process 
$(Y_t=\{L_{2t},L_{2t+1}\})_{i=1}^\infty$ (note the unordered pairs)
is equivalent to the random walk on the complete graph (with self-loops) on vertex set 
$\binom{V}{2}\cup\{\{v,v\}\text{ s.t. }v\in V\}.$
Cover times of this walk will be of particular interest later on, and classical results from the coupon collector problem \cite{erdHos1961classical} will allow us to connect the present work with standard results on matchability under i.i.d.~noise on edges, as discussed in \cite{wu2021settling}; see Section \ref{sec:ER} for details.
We note here that the process $(X_t,C_t)$---where 
$X_t=(L_t,L_{t+1})$ and $C_t$ is the configuration of the lamps in $\{1,\trz\}^{\binom{V}{2}}$ at time $t$---is similar to the classical lamplighter walk (in the language of \cite[Chapter 19]{levin2017markov}) on the line graph of $D_V$, where $D_V$ is the complete digraph on $V$.
The key difference in the present model is that we do not randomize the edge both when visiting and when leaving (as in the classical lamplighter), as well as our chain not being lazy.}
\end{rem}

\section{Mixing and Matching} \label{sec:mix_and_match}

Our aim in this work is to understand how the sequential   error structure introduced by our edgelighter walk process influences graph matchability.
Toward this end, consider the problem of aligning the graph $G_t$ at time step $t$ of the edgelighter walk to the original graph $G_0$.
In particular, our aim is to understand how many steps the walker can proceed before matching fails to recover a significant fraction of the correct vertex alignment.
To formalize the anonymization time we consider the following: 
for $\beta\in(0,1)$, we define the $\beta$-anonymization time, $t_{\text{anon}}^\beta$, via
\begin{equation*}
t_{\text{anon}}^\beta=\min
\left\{t\in\mathbb{Z}^{\geq 0}\text{ s.t. }\underset{P\in\Pi_n}{\text{arg min}}\|G_{t}P-PG_0\|_F\subset \bigcup_{k>n^{\beta}}\Pi_{n,k}\right\},
\end{equation*}
where $\Pi_{n,k}$ is the set of permutations that shuffle exactly $k$ vertex labels for $0<k\leq n$ (note that this is also formalized below in the definition of $\beta$-anonymization in Definition \ref{def:Match}). 
Permutations in the set $\Pi_{n,k}$ are called partial derangements in the combinatorics literature \cite{Stanley2011}.
While we will make this more precise shortly, the intuition is that at such a point, $t_{\text{anon}}^\beta$, much of the structure present in graph $G_0$ has been effectively erased by the noise process. 
At the same time, we seek to tie this anonymization time across the noisy graphs to the \emph{mixing time} of the underlying Markov chain $(Z_t)$ defined in the previous section.

Given an ergodic Markov chain, a branch of modern Markov chain theory is devoted to analyzing how fast  convergence to stationarity occurs, motivating the definition of mixing times of Markov chains \cite{aldous2002reversible,levin2017markov}.
\begin{Def} \label{def:mix}
    Let $(X_n)$ be a discrete time ergodic Markov chain on state space $S$ with one-step transition probability function
    $P : S \times S \rightarrow [0,1]$ and stationary distribution $\pi$.
    The {\em total variation mixing time} of the chain is given by
    \begin{equation*}
    t_m=\min\left\{t\geq 0: \max_{x\in S} \left\|P^t(x,\cdot)-\pi\right\|_{TV}<\frac{1}{4}\right\}.
    \end{equation*}
\end{Def}
\noindent Intuitively, the mixing time captures the number of steps required for the distribution of the Markov chain to become close to its stationary distribution.
We note that the choice of $1/4$ is by convention---it can be replaced with any appropriately small constant; indeed, if we define $t_m(\epsilon)=\min\left\{t\geq 0: \max_{x\in S} \left\|P^t(x,\cdot)-\pi\right\|_{TV}<\epsilon\right\}$, then $t_m(\epsilon)\leq \lceil \log\epsilon^{-1}\rceil t_m$. 
Mixing time analysis arises in a range of applications, including statistical physics, randomized algorithms and probabilistic combinatorics \cite{bayer1992trailing,wilson2003mixing,morris2014mixing,levin2010glauber,banisch2016voter,levin2017markov}.
As a result, a variety of methods have been developed to analyze mixing times, such as spectral gap analysis and coupling methods \cite{aldous2002reversible}. 

Mixing time analysis has also become an important task in the context of analyzing random walks on random networks.
For example, \cite{sood2004first} establishes mixing times for time-homogeneous Markov chains on Erd\H{o}s-R\'{e}nyi graphs, whereas \cite{jafarizadeh2020fastest} analyzes how specific graph structures can be leveraged to achieve faster mixing.
In this light, we will consider the mixing time analysis of the edgelighter graphs on Erd\H{o}s-R\'{e}nyi and SBM random graphs (i.e., $G_0$ is drawn from an Erd\H{o}s-R\'{e}nyi or SBM model).
We surmise here that at a time greater than the mixing time of the Markov Chain $(Z_t)$, the edgelighter has randomized much of the original edge signal in $G_0$.
Indeed, we will make precise shortly the intuition that after edgelighter mixing, the edge probabilities have sufficiently collapsed towards stationarity which results in matchability loss.

To capture this notion of matchability loss, we consider the following definitions, which aim to characterize our ability to recover the vertex labels of a shuffled $G_t$ via matching it to $G_0$.
To achieve this, we seek conditions on $t$ such that we can still correctly match $G_t$ with $G_0$. 
Asking that a ``perfect matching'' that recovers the true correspondences of all $n$ vertices will be too stringent, and thus we here deploy a looser notion of matchability.

\begin{Def}[Ratio Preserved Matchings] \label{def:Match}
Let $G_1 = (V_1,E_1)$ and $G_2 = (V_2,E_2)$ be two vertex-aligned graphs, with respective adjacency matrices $A$ and $B$.  
Let $\alpha,\beta \in (0,1)$ be constant.
We then define
\begin{enumerate}
\item {\bf $\alpha$-matchability.} 
We say that $G_1$ and $G_2$ are $\alpha$-matchable if 
\begin{equation*}
\Prb \left[ \exists P\in \bigcup_{k\geq n^\alpha}\Pi_{n,k} \;\;\text{s.t.}\;\, \tr{APBP^T} \geq \tr{AB} \right] = o(1).
\end{equation*}
In other words, we say graphs $G_1$ and $G_2$ are $\alpha$-matchable if the probability of finding a permutation that shuffles at least $n^\alpha$ nodes and provides at least as good a matching as the true identity mapping as measured by the GMP objective function is asymptotically vanishing.
\item \textbf{$\beta$-anonymization.} We say that $G_1$ and $G_2$ are $\beta$-anonymized if
        \begin{equation*}
        \Prb \left[ \argmax_{P\in\Pi_n} \tr{APBP^T} \subset \bigcup_{k=0}^{n^\beta}\Pi_{n,k}\right] = o(1).
        \end{equation*}
        In other words, graphs $G_1$ and $G_2$ are $\beta$-anonymized if the probability that the optimal matching between them (as measured by the GMP objective function) shuffles at most $n^\beta$ nodes is asymptotically vanishing.
        In other words, with high probability, the optimal matching will shuffle more than $n^\beta$ nodes.
    \end{enumerate}
\end{Def}

\section{Mixing and matching in the Erd\H{o}s-R\'{e}nyi setting}
\label{sec:Main}

We aim to explore the relationship between the mixing time of the edgelighter Markov chain described above and the matchability of the initial graph $G_0$ and the graph $G_t$ at time $t$; recall that $A_0$ and $A_t$ denote their corresponding adjacency matrices, respectively.
We consider first the standard edgelighter walk of Definition \ref{def:stdwalk}.

Note that letting $L_t$ represent the vertex visited by the lamplighter at time $t$, the sequence $\{Y_t=\{L_{2t},L_{2t+1}\}\}_{t=0}^\infty$ forms a Markov Chain on state space $S = \binom{V}{2}\cup\{\{v,v\}:v\in V\}\}$; see Remark~\ref{rem:global}. 
Let $\tcov$ be the cover time of the Markov Chain $(Y_t)$, i.e., 
\begin{align}
\label{eq:cov}
\tau_{\text{cover}}=\min\{t:S\subset \cup_{\ell\leq t} Y_\ell\}; \quad \tcov=\max_{s\in S}\mathbb{E}_s(\tau_{\text{cover}}). 
\end{align}
where $\mathbb{E}_s(\cdot)$ denotes the expectation where the lamplighter walk starts in state $s$.
Note that simple coupon collector asymptotics give us that the cover time for this chain, $\tcov$, is of order $\Theta(n^2\log n)$, and that $\tau_{\text{cover}}$ is of order $\Theta(n^2\log n)$ with high probability.
If $q_1$ and $q_2$ are chosen so that $q_1=1-q_2$, then for $t \ge \tau_{\text{cover}}$, we have that $G_t$ is distributed as $\ER(n,q_2)$ and that $G_t$ is independent of $G_0$.  
In this case, it is clear that 
$$\mathbb{P}\left(I_n\in \argmin_{P\in\Pi_n}\|A_0-PA_tP^T\|_F\,\Big|\,t\geq \tau_{\text{cover}}\right)$$ is bounded away from 1 for all $t>\tau_{\text{cover}}$, as by this time any permutation is equally likely to be optimal.

Our first lemma then relates the mixing time of $(Z_t)$, which we will denote $t^\circ_m$, to a suitably defined cover time of $(X_t)$ where $X_t=(L_t,L_{t+1})$.
Define 
$\tau^{(u)}_{\text{cover}}$ via
\begin{equation*}
\tau^{(u)}_{\text{cover}}=\min\left\{t:\binom{V}{2}\subset \bigcup_{\ell\leq t-1}\{L_\ell,L_{\ell+1}\}\right\},
\end{equation*}
and define $t^{(u)}_{\text{cover}}=\max_x \mathbb{E}_x \tau^{(u)}_{\text{cover}}$.
Note that it is immediate that $\tau^{(u)}_{\text{cover}}\leq \tau_{\text{cover}}$, where $\tau_{\text{cover}}$ is as defined in Eq. \ref{eq:cov}; note that the proof of the following Lemma can be found in Appendix \ref{app:covermix}.
\begin{lem}
\label{lem:covermix}
With notation as above, let $q_1=1-q_2$ in the edgelighter walk of Definition \ref{def:stdwalk}.
We then have that $\frac{1}{8}t^{(u)}_{\text{cover}}\leq t_m^\circ\leq 11t^{(u)}_{\text{cover}}$.
\end{lem}
\noindent Simple coupon collector asymptotics (see Appendix \ref{app:pijt} for proof) give us that for sufficiently large $n$,
\begin{equation}
\label{eq:covn2}
\frac{1}{2}n^2\log n\leq t_{\text{cover}}^{(u)}\leq \frac{5}{2} n^2\log n,
\end{equation}
and thus $t_m^\circ=\Theta(n^2\log n)$, where we recall $t_m^\circ$ is the mixing time of $Z_t$. 

Our focus now shifts to exploring models for $G_0$ and the edgelighter walk where signal recovery is possible up to time $\tcov$ and identifying cases where signal recovery fails before time $\tcov$.
We begin by considering the case where $G_0$ is distributed as an ER random graph.

\subsection{Standard Edgelighter Walk on Erd\H{o}s–R\'{e}nyi Graphs}
\label{sec:ER}

We first consider the Standard edgelighter walk on an ER graph $G_0 = (V,E)\sim \ER(n, p)$.
We will assume that the initial starting vertex for the edgelighter is chosen uniformly at random from $V$.
Letting $G_t$ be the graph after $t$ steps of the edgelighter walk, we consider here the parameters $q_1=(1-p)$, $q_2=p$ so that $G_t\sim \ER(n, p)$ as well.

To compute the edge-level correlation between $G_t$ (with adjacency matrix $A_t$) and $G_0$ 
(with adjacency matrix $A_0$) in this setting, define
\begin{equation*}
\mathfrak{p}_{t,u,v}
=\Prb\bigg(\{u,v\}\notin \cup_{s\leq t-1}\{L_s,L_{s+1}\}\bigg);
\end{equation*}
i.e., $\mathfrak{p}_{t,u,v}$ is the probability that the edge $\{u,v\}$ has not been traversed by the edgelighter by time $t$.
In Appendix \ref{app:pijt}, we show that
    \begin{equation}
\label{eq:pijt}
    \mathfrak{p}_{t,u,v} = \exp\left\{ -\Theta\left(\frac{t}{n^2}\right)\right\}
\end{equation}
After the edge $\{u,v\}$ has been traversed, its on/off status is independent of the initial on/off state, and we thus have
\begin{equation*} \begin{aligned}
\Corr(A_{0,u,v},A_{t,u,v})
&=\frac{\mathbb{E}(A_{0,u,v}A_{t,u,v})-p^2 }{ p(1-p)} =\frac{p\mathfrak{p}_{t,u,v}+p^2(1-\mathfrak{p}_{t,u,v})-p^2 }{ p(1-p)} =\mathfrak{p}_{t,u,v}
\end{aligned} \end{equation*}
We note, however, that unlike the traditional correlated Erd\H os-R\'enyi setting where edges between different vertex pairs are independent across graphs, the edgelighter introduces a nuanced dependency structure across the networks $G_0$ and $G_t$. 
For example, at time $t$, if there are $t$ edges in $G_t$ with opposite on/off status to those in $G_0$, then the remaining edges in $G_t$ must be identical to those in $G_0$. 
Nonetheless, we suspect that this dependence is mitigated in time, and conjecture that a matchability phase transition occurs at edge-wise correlation of order $\sqrt{n^{-1} \log n}$ (the same order as in the correlated Erd\H os-R\'enyi setting \cite{wu2021settling,lyzinski2016information,lyzinski2014seeded}), which would imply that $t=\Omega(n^2\log n)$ is sufficient for the graphs to be effectively anonymized and that if $t=O(n^2\log n)$, matchability is preserved.

With the matchings defined in Definition~\ref{def:Match}, we next provide a sequence of results that give (at least a partial) affirmative answer towards the matchability thresholds conjectured above.
In the following theorems, we consider $G_0\sim ER(n,p)$ and consider $1-q_1=q_2=p$.
At each step, we again match $G_t$ with $G_0$, and we have the following theorem (proven in Appendix \ref{app:thm1}).
\begin{thm} \label{thm:global}
Let $\{G_t\}_{t\in\mathbb{N}}$ be as described above with $1-q_1=q_2=p$. 
Let $\alpha,c>0$ be fixed constants such that $\alpha>5c$. 
If $t\leq c n^2\log n$, then $G_0$ and $G_t$ are $\alpha$-matchable; specifically,
\begin{equation*}
\PP\left[ \exists\, P\in \bigcup_{k\geq n^\alpha}\Pi_{n,k}\text{ s.t. }\tr{A_0PA_tP^T}\geq \tr{A_0A_t}\right]
\leq  2\exp\left\{ -\Omega\left(\frac{n^{2\alpha-5c} }{\log n}\right) \right\} .
\end{equation*}
\end{thm} 
\noindent 
With notation as above, for the matching anonymization bound we consider the time $\tau^{(u)}_{\text{cover}}$, as after this time, $G_{t}$ is effectively $\text{ER}(n,p)$ independent of $G_0$ and any permutation is equally likely to be optimal.
As noted in Remark \ref{rem:global}, the (unordered) pairs $(\{L_{2i},L_{2i+1}\})_{i=1}^\infty$ are independent and uniformly distributed over $S=\binom{V}{2}\cup\{\{v,v\}:v\in V\}$.
With $\tau_{\text{cover}}$ as in Eq. \ref{eq:cov}, standard coupon collector asymptotics yield that (as $|S|\leq n^2$)
\begin{equation} \label{eq:couponcollector}
\PP\left[ \tau^{(u)}_{\text{cover}}\geq 4 n^2\log(n) \right]
\leq \PP\left[ \tau_{\text{cover}} \geq 4 n^2\log(n) \right]
\leq n^2 e^{-4\log n}=n^{-2}.
\end{equation}
We then have the following theorem.
A proof is given in Appendix \ref{sec:thm2}.

\begin{thm} \label{thm:globalbad}
Let $\{G_t\}_{t\in\mathbb{N}}$ be as described as in Theorem \ref{thm:global} 
with $1-q_1=q_2=p$. 
If $t\geq 4 n^2\log n$, then for any $\beta\in(0,1)$, $G_0$ and $G_t$ are $\beta$-anonymized; specifically,
\begin{equation*}
\PP\left[ \argmax_{Q\in\Pi_n} \tr{A_0QA_tQ^T}\subset \bigcup_{k=0}^{n^\beta}\Pi_{n,k}  \right]\leq  n^{-2}+(2n^\beta)^{1-n/(2n^\beta+1)}
\end{equation*}
\end{thm}
\noindent Theorems \ref{thm:global} and \ref{thm:globalbad} 
establish that for the standard edgelighter, de-anonymization between $A_t$ and $A_0$ occurs on the same order as the mixing of $(Z_t)$, namely in order $t=\Theta(n^2\log n)$ steps.  
\begin{rem}
\emph{We note that Theorems \ref{thm:global} and \ref{thm:globalbad} 
hold (after appropriately adjusting constants) if $1-q_1=q_2=q$ for any $q\in(0,1)$. 
We use $q=p$ for consistency with the $\ER(n,p)$ model parameters before and after $\tau^{(u)}_{\text{cover}}$.}
\end{rem}

\section{(Partial) Anonymization Before Mixing for Structured Graphs}
\label{sec:SBM}

We have observed that for $G_0$ sampled from a simple ER model, there is an agreement between the mixing time of the noise process $(Z_t)$ and the anonymization time.
Herein, we seek to identify models where the matching fails (i.e., the graph has been sufficiently anonymized as in Definition \ref{def:Match}) before the noise mixes the distribution of the network.
To explore this, we consider a more structured distribution on graphs: the Stochastic Block Model (SBM).
The specific model we will examine is 
\begin{equation*}
\SBM\left(K_n , \Lambda_n=[\Lambda_{i,j;n}], \tau_n \right)\subseteq \mathcal{G}_n,
\end{equation*}
Within this SBM model, we define the edgelighter walk as follows
\begin{Def}[(Block Edgelighter Walk)]
\label{def:sbmwalk}
Let $G_t = (V, E_t)$ be the state of the graph at time $t$, and assume that the vertex set of $G_0$ is partitioned into $K$ communities $B_1,B_2\dots,B_K \subseteq V$ (e.g., $G_0$ could be sampled from the SBM model specified above). 
The edgelighter here is again modeled by a time homogeneous Markov chain $(W_t)$, this time on state space $\mathcal{C}\times V\times \{1, \trz\}^{\binom{V}{2}}$, where $\mathcal{C} = [K]$ is the collection of all community labels in the graph. 
If the lighter is at community $i$ (i.e., is in $B_i$) and vertex $u\in B_i$, then for the next step, the edgelighter stays or leaves community $B_i$ with equal probability (i.e., equal to 1/2).
\begin{itemize}
\item[i.] If the edgelighter selects to stay in community $B_i$, then the lighter selects a vertex $ v\in B_i$ uniformly at random, and moves to $v$. 
If $v=u$, then $G_{t+1}=G_t$ and nothing is done.
If the edgelighter moves to vertex $v\neq u$ at time $t+1$, then we consider two cases.
If $h_{G_{t}}(\{u,v\})=1$, then the graph evolves via 
$G_{t+1}=(V,E_{t+1})$ where 
\begin{align}
h_{G_{t+1}}(e)=\begin{cases}
h_{G_{t}}(e)& \text{ if }e\neq \{u,v\};\\
X_t(1-h_{G_{t}}(e))+(1-X_t)h_{G_{t}}(e)& \text{ if }e= \{u,v\}
\end{cases}
\end{align}
        where $X_t\sim\Bern(q_{i;1})$ is  independent the edgelighter walk and of $(G_s)_{s\leq t}$.
             If $h_{G_{t}}(\{u,v\})=\trz$, then the graph evolves via 
        $G_{t+1}=(V,E_{t+1})$ where 
        \begin{align}
        h_{G_{t+1}}(e)=\begin{cases}
        h_{G_{t}}(e)& \text{ if }e\neq \{u,v\};\\
            Y_t(1+h_{G_{t}}(e))+(1-Y_t)h_{G_{t}}(e)& \text{ if }e= \{u,v\}
        \end{cases}
        \end{align}
        where $Y_t\sim\Bern(q_{i;2})$ is  independent the edgelighter walk and of $(G_s)_{s\leq t}$.   
        \item[ii.] If the edgelighter selects to leave $B_i$, then the lighter first choose a new community $B_j\neq B_i$ uniformly at random from all $K_n-1$ non-$B_i$-communities.
        The edgelighter then independently selects an edge $e_1$ between $B_i$ and $B_j$ uniformly at random and a non-edge $e_2$ between $B_i$ and $B_j$ uniformly at random, and sets $h_{G_{t+1}}(e_1)=\trz$ and $h_{G_{t+1}}(e_2)=1$.
        The edgelighter then pick a node $w\in B_j$ uniformly at random from those in $B_j$, and moves to $w$ at time $t+1$. 
        \end{itemize}
    {We can envision this edgelighter as performing a hierarchical walk on the graph with communities.  At the community level, the edgelighter is moving between communities and randomizing edges connecting the chosen communities uniformly at random.  Within communities, the edgelighter is behaving like the standard edgelighter walk.}
\end{Def}

As in the simple edgelighter chain, the block edgelighter chain is aperiodic, irreducible (on the state space restricted so that the number of edges between any two communities is fixed to be equal to the value in $G_0$), and reversible.
Here the stationary distribution $\pi^\bullet$ is defined via 
\begin{equation*}
\pi^\bullet(B_i,u,\vec c)\propto \frac{1}{K} \frac{1}{|B_i|}\prod_{j}q_{j;2}^{m_{j,\vec{c}}}q_{j;1}^{r_{j,\vec{c}}} \frac{1}{\prod_{i,j}\binom{n_in_j}{m_{ij}}},
\end{equation*}
where $K$ is the number of communities, $u\in B_i$, and $\vec c$ is the edge configuration with $m_{i,\vec{c}}$ the number of $1$'s present between vertices in $B_i$ and $r_{i,\vec{c}}$ the number of $\trz$'s present between vertices in $B_i$, $n_1$ the number of vertices in $B_i$ and $m_{ij}$ the number of edges between $B_i$ and $B_j$; for a proof of this form of the stationary distribution, see Appendix \ref{app:stat2}.
With the Block Edgelighter Walk defined, we have the following result on \emph{local $\beta$-anonymization of a community in the graph} before the block edgelighter noise model has fully mixed.  
We say that a community $B_i$ of size $n_i$ has been $\beta$-anonymized at time $t$
if the following holds (where $\Pi_{B_i,\beta}$ is the set of permutations in which at most $n_i^\beta$ of the vertices in $B_i$ have been shuffled)
        \begin{equation*}
        \Prb \left( \underset{P\in\Pi_n}{\text{arg max}}\operatorname{Tr}(A_0PA_tP^T) \subset \Pi_{B_i,\beta}\right) = o(1).
        \end{equation*}
Note that the proof of Theorem \ref{thm:SBM} can be found in Appendix \ref{sec:pf3}. 
\begin{thm} \label{thm:SBM}
Assume there exist constants $b_1,b_2,b_3,b_4>0$ and $0<a_1< a_2<1$ such that, assuming without loss of generality that $n_1=\min n_i$,
\begin{align*}
    n_1=b_1n^{a_1};&\quad 
    n^*:=\max_i n_i=b_2n;\quad \Lambda_{1,j;n}\leq b_3\frac{n^{a_2}}{n_jn_1}\text{ for all }j\neq 1;\quad K_n=b_4;\\
    \Lambda_{i,j;n}&\geq 5\frac{\log n}{n_in_j}\text{ for all }i,j\text{ with }i\neq j;\quad \Lambda_{i,i;n}\geq 5\frac{\log n}{\binom{n_i}{2}}\text{ for all }i
    \end{align*}
Consider the SBM lamplighter defined above with $1-q_{i;1}=q_{i;2}=q$ for all $i$.    
Define the event
\begin{equation*}
\mathcal{E}_{t,\beta}=\left\{\underset{P\in\Pi_n}{\text{arg max}} \operatorname{Tr}(A_0PA_tP^T) \subset \Pi_{B_1,\beta}\right\},
\end{equation*}
recalling that $\Pi_{B_1,\beta}$ is the set of permutation matrices that shuffles at least a fraction of $\beta$ labels in community $B_1$,
    and letting $N = n^{1+a_2}\log^{3/2}(n)$.
We then have for a suitable constant $C>0$, 
\begin{align*}
    \mathbb{P}(\mathcal{E}_{N,\beta})&
    \leq O\left((2n_1^\beta)^{-n_1/(2n_1^\beta+1)}\right)+O\left((n_1n^*)^{-Cn^{a_2}}\right)+O( n^{2a_2}\log\left(n\right)b_4^{-n^{(1-a_2)/2}})\\
  &\quad+O(n^{a_2-1}/2))+e^{-\Omega(n^{a_1}\log ^{1/2}(n))}
\end{align*}
\end{thm}

\noindent We shall refer to the stage where at least one of the communities of the network are $\beta$-anonymized (community 1 in our setup), as the point when the network is locally anonymized. 
Theorem \ref{thm:SBM} combined with Lemma \ref{lem:covermix2} below shows that, 
under mild assumptions on the block edgelighter on an SBM (the key being that there is a community that is of lower order in size compared to the largest community), the matching corrupts locally before the Markovian noise globally mixes.
We note that Theorem \ref{thm:global} does not apply here because the community structure introduces different probabilities for walking within or between communities.
The following Lemma, proven in Appendix \ref{app:covermix2}, completes this argument:
\begin{lem}
\label{lem:covermix2}
With assumptions as in Theorem \ref{thm:SBM}, let $t^\bullet_m$ be the mixing time of the Block Edgelighter Walk.
We then have that $t^\bullet_m=\Omega(n^2)$.
\end{lem} 
\noindent The mixing time $t^\bullet_m=\Omega(n^2)$, and hence after 
$n^{1+a_2}\log^{3/2}(n)$ steps (sufficient for local $\beta$-anonymization), the chain $(W_t)$ has not globally mixed.

The local $\beta$-anonymization proven in Theorem \ref{thm:SBM} is a consequence of \emph{locally mixing} the smallest community in a sense; indeed, at the time of local $\beta$-anonymization, the edges within and across this smallest community have been effectively mixed (i.e., are close to the appropriate stationary probabilities restricted to these edges).  
Lastly, although for the Block Edgelighter Walk, matching locally corrupts earlier than the mixing time of the noise, Theorem \ref{thm:SBM} only provides a condition for the inability to achieve exact matching as defined in \cite{wu2021settling}. According to Theorem \ref{thm:SBM}, after $\Theta(n^{1+a_2}\log(n))$ steps, we are only guaranteed anonymization of a constant number of communities (at least 1), not the entire network. We hypothesize that ensuring global $\beta$-anonymization still requires $\Omega(n^2\log(n))$ time, and further analysis of such bounds is a topic of future work. 

\begin{rem}
        \emph{In Theorem \ref{thm:SBM}, the assumptions on the parameters in the SBM can be tightened.
        In essence, the assumptions need to ensure that the smallest community $n_1$ is of lower order than the largest community, and there are not too many edges connecting the smallest community to the rest of the graph.
        As we are interested in this example from a heuristic perspective, we do not pursue necessary and sufficient conditions herein.}
\end{rem}

\section{Experiments}
\label{sec:experiments}

To explore our theories further, we present two sets of experiments.
In the first, the edgelighter is walking on synthetic random base graphs (an Erd\H os-R\'enyi or SBM base graph), and in the second the edgelighter is walking on real data base networks $G_0$. 
In these simulations, we approximate both a suitable anonymization time and the edgelighter walk's mixing time, with the aim of showing approximate anonymization occurring at the time of mixing.
Note here that as the precise mixing time is hard to quantify, we employ Lemmas \ref{lem:covermix} and \ref{lem:covermix2} in the standard edgelighter walk to equate (at least asymptotically) the mixing time $t_m^\circ$ and the cover time $\tcov$ of the walk; our amended aim is then showing approximate anonymization occurring at the time of covering.
In the SBM/block setting, the cover time provides a lower bound on mixing. 

To address the computational intractability of solving the exact graph matching problem, we use the Seeded Graph Matching (SGM) algorithm introduced by \cite{ModFAQ} as follows.
To expedite the experiments on matchability for the SBM walk and real data, we initialize the matching algorithm at the ground truth---the identity matrix. 
Note that when anonymization occurs, the algorithm's optimizer should move away from the ground truth. 
While the ground truth being a local optimum can prevent this move, this initialization provides a good heuristic for assessing the optimality of the ground truth, while effectively saving computational resources.
The SGM algorithm \cite{ModFAQ, shirani2017seeded} further makes use of seeded vertices (i.e., vertices whose correspondence is known a priori and fixed throughout the matching process) to enhance matching performance.
In our experiments, our seeds are randomly selected from the set of nodes, and the algorithm ensures that the matching results map these seeded nodes to themselves. Unless stated otherwise, we select 5\% of the nodes, uniformly at random, as seeds.

\subsection{Simulated edgelighter walks on ER and SBM graphs}
\label{exp:ER}
\begin{figure}[t!]
    \centering
    \includegraphics[width = \textwidth]{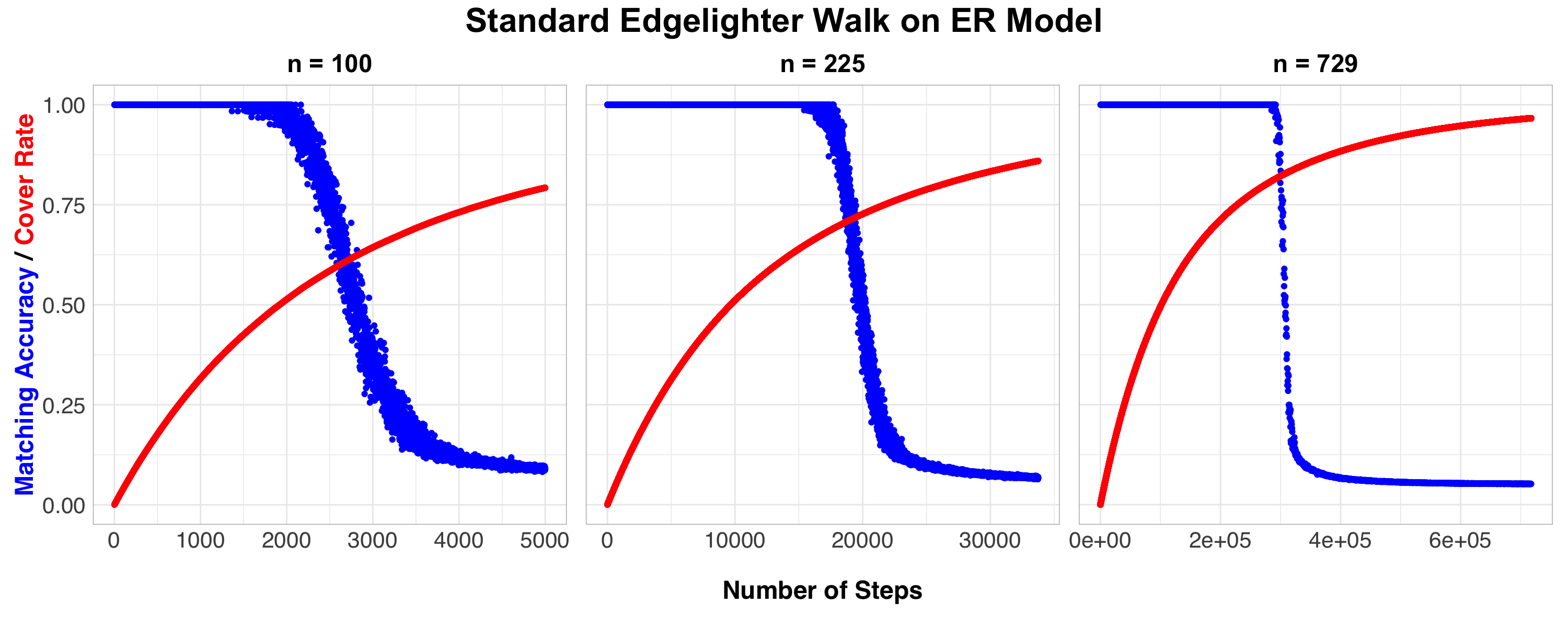}
    \includegraphics[width = 0.88\textwidth]{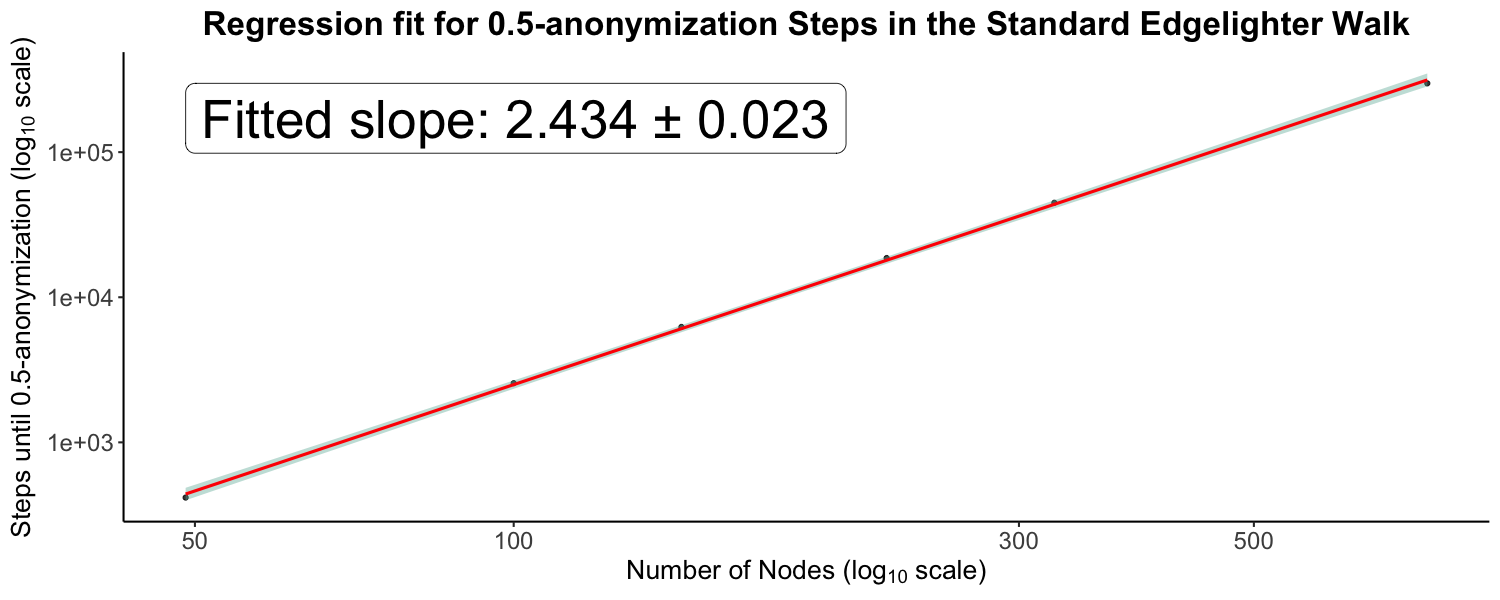}

    \caption{On the top panel, we plot matching correctness vs. number of steps for $n=100$ (left), $n=225$ (middle) and $n=729$ (right) in blue. For all these plots, we further plot (in red) the cover rate of the edges vs number of steps. In the bottom panel, we plot a regression of number of steps needed until a $0.5$-anonymization happens as a function of the number of nodes $n$ on a log-log scale.}
    \label{fig:global_ratio}
\end{figure}

We begin by examining the standard edgelighter walk model of Definition \ref{def:stdwalk} on ER graphs. 
We sample ER initial graphs with $p=0.5$ and node counts $n=49,\, 100,\, 144,\, 225,\, 324,$ and $729$; here $q_1=q_2=1/2$. 
For each graph, we perform a standard edgelighter walk. For smaller graphs ($n=49,\,100,\,144$), we check matching correctness after each move. To save computational resources for larger graphs ($n=225,\, 324,\, 729$), we check matching after a fixed number of steps, specifically after $s_n = 3,\, 30,$ and $300$ steps for $n=225,\, 324,$ and $729$ respectively. 
In the top panel of Figure \ref{fig:global_ratio}, we plot (in blue) the matching correctness versus the number of steps for $n=100,\,225$ and $n=729$. Results for the other values of 
node counts can be found in Appendix \ref{app:ER_plt}. Additionally, in red we also plot the cover rate of the edges against the number of steps for these experiments.
In the bottom panel of Figure \ref{fig:global_ratio}, we show a log-log plot of the number of steps needed until $0.5$-anonymization happens vs the number of nodes, along with a fitted regression line. 
Additional plots in log-log scale for the number of steps needed until the $0.25$-anonymization and until the $0.75$-anonymization for the same run can be found in Appendix \ref{app:ER_plt}.
The slope of the fitted line is slightly over 2, which aligns with our theoretical result that $\Theta(n^2\log n)$ steps are required to corrupt the matchability.

Next, we consider the more structured block edgelighter on $G_0$ sampled from a Stochastic Block Model (SBM) to illustrate a scenario where anonymization occurs before mixing. 
Here, we sample graphs from the following SBM model: 
\begin{equation*}
    \SBM\left(K \!\!=\!\! \left\lfloor\! \frac{n-n^{1/4}-n^{4/3}}{n^{2/3}}\!\right\rfloor\!, ~ \Lambda \!=\! \frac{\log(n)}{n^{3/4}}J + \operatorname{diag}(1/2), ~\tau \!= \!\tau_n \right)
\end{equation*}
where $\tau_n$ is encoded by a vector of block sizes $\vec n =(n_1,\cdots,n_K)= (n^{1/4}, n^{2/3}\cdot \mathbf{1}_{K-2}^T, n^{3/4})$ according to $\tau_n(j) = \sum_{k=0}^{K-1}\mathbb{1}_{\left\{j\geq \sum_{i=1}^{k}{n}_i\right\}}$ where $ n_0=0$. 
These parameters were chosen to reflect the key conditions of Theorem \ref{thm:SBM}, namely the community imbalance, and sparsity of connection between communities 1 and the remainder of the graph; note that an analogue of Theorem~\ref{thm:SBM} can be proven under these conditions.
For the experiments, we select node counts $n=81,\, 256,$ and $625$ and perform the SBM edgelighter walk on these sampled networks. Matching is conducted after every step for $n=81$, after every 90 steps for $n=256$, and after every 2100 steps for $n=625$. Additionally, recall that we initialize the matching algorithm at the identity matrix.

\begin{figure}[t!]
    \centering
    \includegraphics[width = \textwidth]{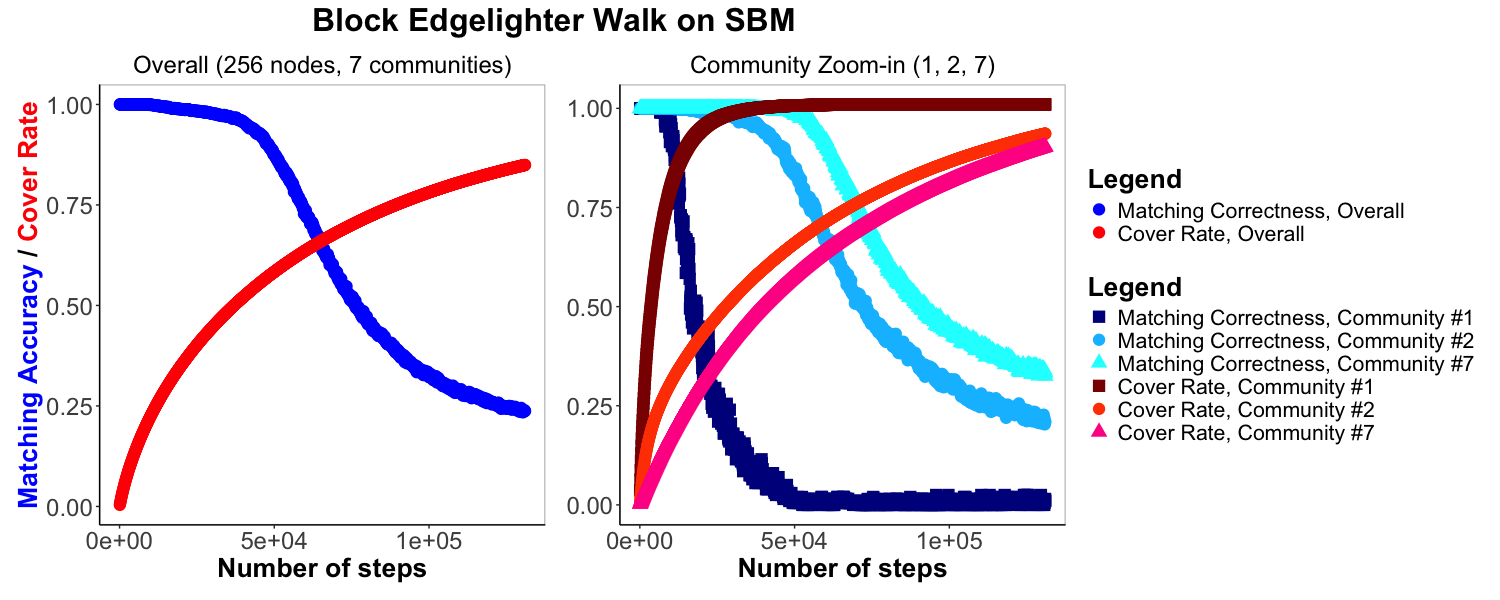}
    \includegraphics[width=0.88\linewidth]{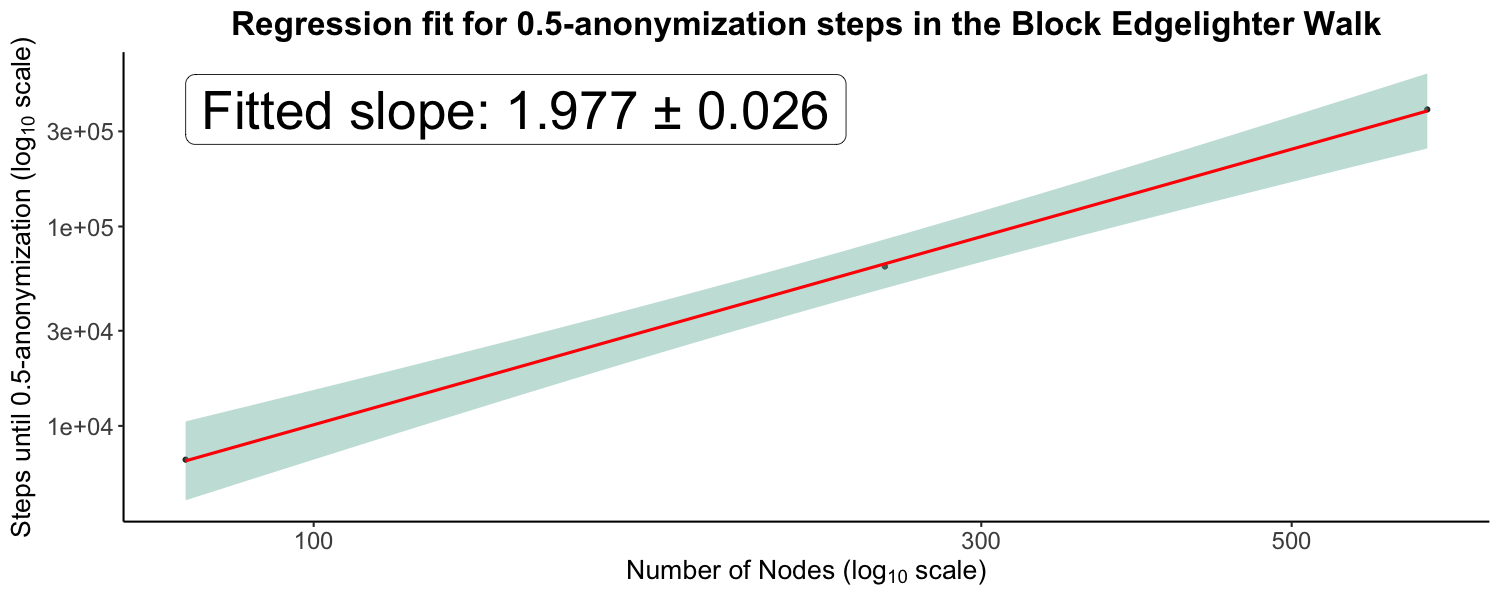}
    \caption{We show the matching correctness as a function of the number of steps for $n=256$ vertices under an SBM-distributed initial graph under a Block Edgewalker Model.
    The top-left panel shows the matching correctness (in blue) versus the number of steps for the entire network; the top-right panel shows matching correctness versus the number of steps for community 1 (the smallest community), community 2 (a randomly chosen community from among the $K-2$ communities of size $n^{2/3}$, here $K=7$) and community $7$ (the largest community). In all these plots, we also include the edge cover rate (in red) against the number of steps. In the bottom panel, we plot the log-log graph of the number of steps required until a $0.5$-anonymization is achieved versus the number of nodes, with a fitted regression line superimposed on these points.
    }
    \label{fig:SBM}
\end{figure}

In Figure \ref{fig:SBM}, we present selected plots for the case of $n=256$ nodes; plots for other communities in the $n=256$ case, as well as additional plots for the $n=81$ and $n=625$ cases, can be found in Appendix \ref{app:SBM_plt}. 
In the top panel of Figure \ref{fig:SBM}, we show the matching correctness in blue versus the number of steps for the entire network as well as three selected communities: community 1 (the smallest community); community 2 (a randomly chosen community with size $n^{2/3}$); and community $7$ (the largest community). In all these plots, we also include the edge cover rates in red against the number of steps. 

We observe that the smallest community exhibits noticeably faster anonymization, which aligns closely with our proposed theoretical framework. Additionally, the largest community demonstrates a slower anonymization rate compared to the other communities. It is also important to note that, when comparing this SBM edgelighter setting to the standard edgelighter setting on ER graphs, the transition away from matchability appears to be less sharp.
We also plot in the bottom panel of Figure \ref{fig:SBM} a log-log graph of the number of steps required until a {$0.5$}-anonymization is achieved versus the number of nodes, with a fitted regression line superimposed on these points. 
Additional plots in log-log scale for the number of steps needed until the $0.25$-anonymization and until the $0.75$-anonymization for the same Block Edgelighter run can be found in Appendix \ref{app:SBM_plt}.
As observed, the slope of the regression line is less than 2, and is about 0.5 less than the slope obtained for the standard edgelighter walk on ER networks. This value aligns with the suggested by our theoretical framework in Section \ref{sec:SBM}.

\subsection{Real data experiments}
\label{sec:data}

\begin{figure}[t!]
    \centering
    \includegraphics[width=0.4\linewidth]{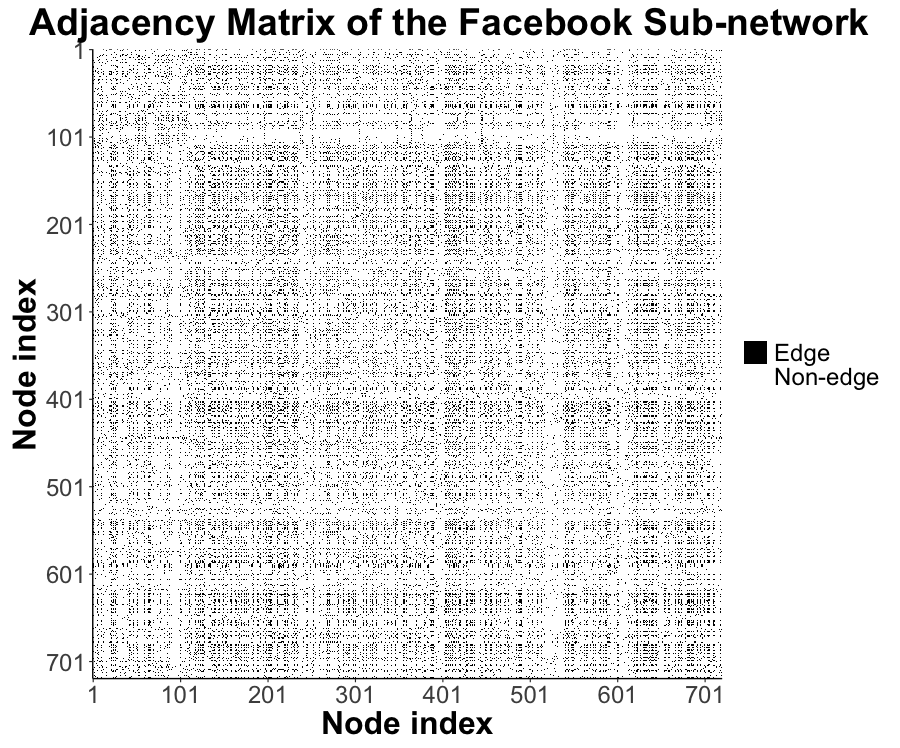}\hspace{3pt}
    \includegraphics[width=0.45\linewidth]{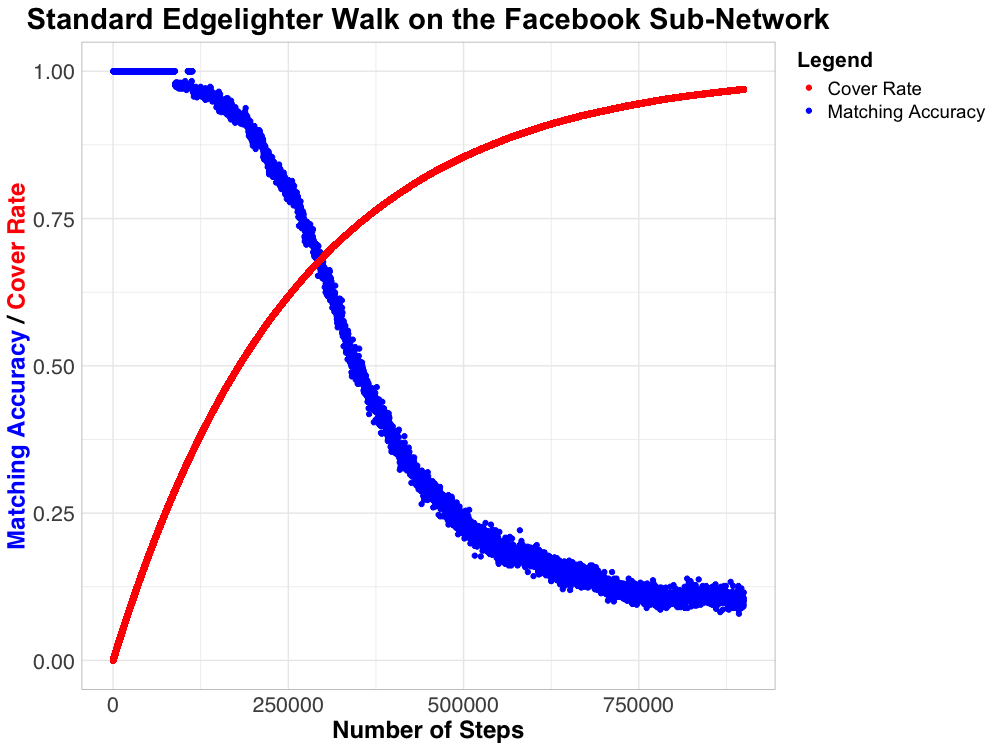}
    \caption{Left: The adjacency matrix of the induced Facebook subgraph that we perform our edgelighter walk on. Right: matching correctness (blue) and cover rate (red) as a function of the number of steps for the Facebook Sub-network. }
    \label{fig:Facebook}
\end{figure}

To validate the proposed theoretical results on a real-world network, we apply the standard edgelighter model to a friendship network from \cite{mcauley2012learning}. 
We downloaded the network from \url{https://snap.stanford.edu/data/ego-Facebook.html}, and extract the induced subgraph containing nodes $1921$ through $2640$ ( $n=720$ nodes in total). A plot of the corresponding adjacency matrix is shown in the left panel of Fig. \ref{fig:Facebook}. 
On the selected network, we run the standard edgelighter walk for 900,000 steps, performing matching to the original network every 150 steps. 
The edgelighter walk introduces noise that can be interpreted as individuals forming new friendships or losing contact with those they have not interacted with recently. We also track the covering rate after each step, defined as
$$\frac{\text{\# distinct edges traversed after the step} }{\binom{720}{2}}.$$
In the right panel of Fig. \ref{fig:Facebook}, we plot the matching accuracy versus number of steps curve in blue, alongside the covering rate versus number of steps curve in red. The results show a similar pattern to what we observed for the standard edgelighter walk on Erd\H{o}s–R\'{e}nyi graphs, as seen in Fig. \ref{fig:global_ratio}.

\begin{figure}[t!]
    \centering
    \includegraphics[width=0.4\linewidth]{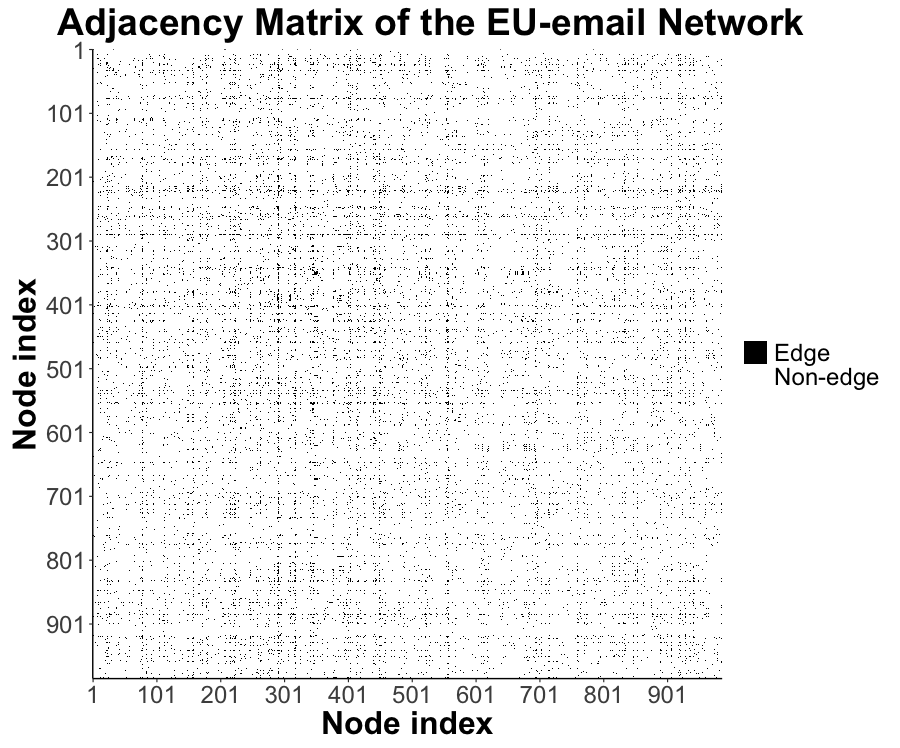}
    \caption{The adjacency matrix of the undirected EU Email communication network. Node indices are reordered based on department memberships. }
    \label{fig:EU_Email_Adj}
\end{figure}

\begin{figure}[t!]
    \centering
    \includegraphics[width = \textwidth]{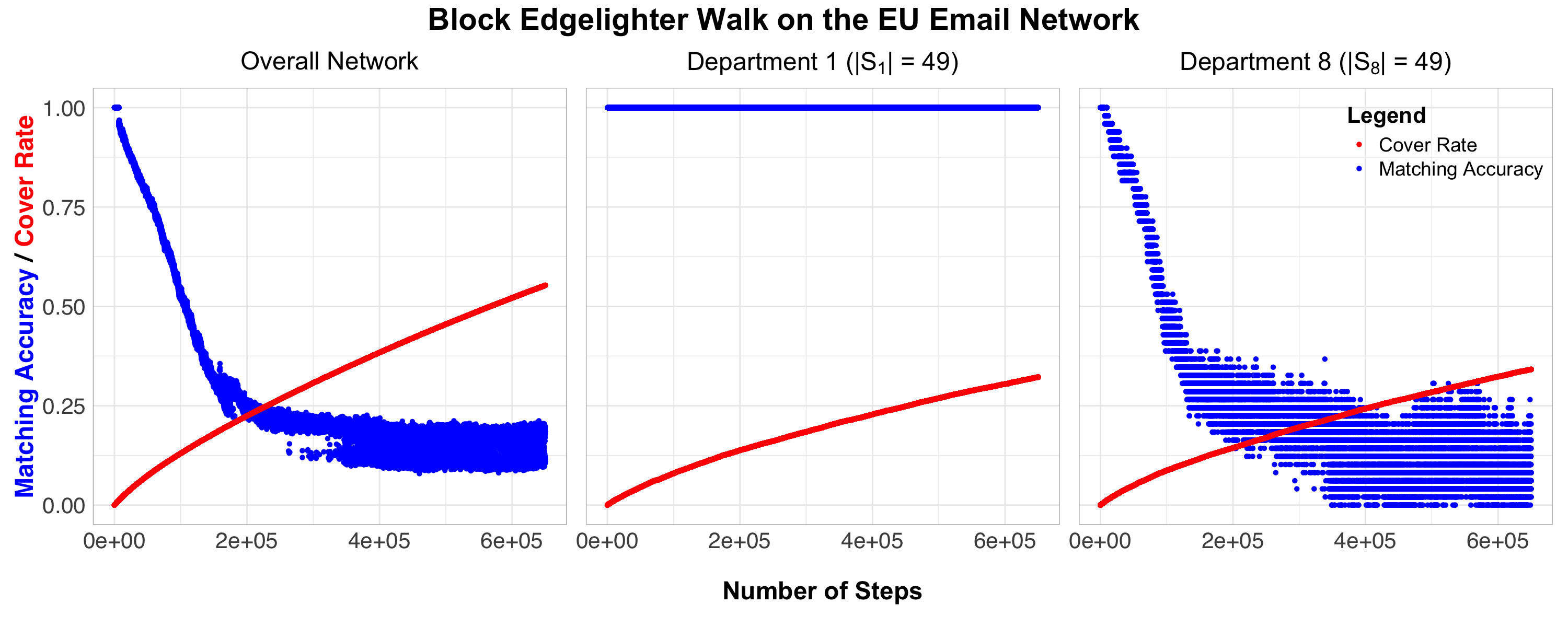}
    \caption{We plot the matching correctness vs iteration plot (in blue) and the cover rate vs iteration curve (in red) for the entire undirected EU Email Network (left panel) and for the sub-networks induced by nodes from community (department) \#1 (middle panel) and by nodes from community (department) \#8 (right panel)}
    \label{fig:EU_Email}
\end{figure}

To further validate the proposed results on a more structured real-world network, we implement the Block edgelighter walk on an email communication network generated from a large European research institution. In this network, nodes represent individuals, and edges indicate email exchanges between them. 
A numerical department label to which each individual belongs is provided as a node attribute, allowing us to treat each department as a distinct community. For more detailed information about this network and its applications, refer to the original source \cite{leskovec2007graph}.
We obtained the network from \url{https://snap.stanford.edu/data/email-Eu-core.html}. For simplicity and better alignment with our theoretical framework, we only consider the largest connected component of the network and made the network undirected by dropping the directedness of edges. 
The resulting adjacency matrix is shown in Fig. \ref{fig:EU_Email_Adj}, where we reordered the nodes based on their communities. As in the block edgelighter on SBM model experiments discussed earlier, we apply the same block edgelighter walk of Definition \ref{def:sbmwalk} to this network, where we treat departments as communities and for each email user, we use their department memberships as the their ground truth community labels. 
To save computational resources, matching with the original network is performed every 220 steps. In the left panel of Fig. \ref{fig:EU_Email}, we plot the matching correctness (blue) versus number of steps curve, along with the edge cover rate (red) versus number of steps curve. The plot shows a similar pattern to what we observed for the block edgelighter walk on the simulated SBM graph in Fig. \ref{fig:SBM}.

To gain a better understanding of anonymization within communities, the middle and right panels of Fig. \ref{fig:EU_Email} show the matching correctness versus number of steps plot and the cover rate versus number of steps curve for the subgraph induced by nodes from community (department) \#1 (middle panel, 49 members) and by nodes from community (department) \#8 (right panel, 49 members). 
Additional plots for some other selected departments can be found in Appendix \ref{app:EU_plt}. From the plots, we observe that anonymization occurs much earlier in some communities. This observation supports our hypothesis that faster anonymization in structured networks is local and is due to the anonymization of certain communities, while we suspect that globally corrupting the entire matching still requires time on the same order as the time needed for the noise to globally mix. 

\section{Conclusion and Discussion}
In this paper, we explored the relationship between the anonymization time of graph signals and the mixing time of Markovian noise on the graph, across various graph models under different lamplighter-like walk schemes, which we refer to as the edgelighter walk model. By examining both theoretical bounds and simulation results, we established that the time required to disrupt the matchability between an original graph and its noisy counterpart aligns closely with the mixing time of the noise, especially in the context of Erd\H{o}s–R\'{e}nyi random graphs. Specifically, our results demonstrate that for the standard edgelighter walks on these graphs, the matching between the original and noisy graphs deteriorates after $\Theta(n^2\log n)$ steps, consistent with the mixing time of the edgelighter. 
However, in some structured graph models, particularly the Stochastic Block Model (SBM) with unbalanced community structures, local anonymization can occur before the noise has fully mixed. This highlights distinct dynamics in such models. Our findings suggest that in SBM networks with $n$ nodes and communities with sufficiently different sizes, the time required to break matchability locally is approximately $\Theta(n^\alpha\log n)$ for some $\alpha < 2$, significantly earlier than the $\Theta(n^2\log n)$ threshold for dense graphs.
Note that the mixing time depends on the number of edges in the graph. Thus, for sparser graphs (with fewer than $\Theta(n^2)$ edges), the mixing time of the noise could be smaller than $\Theta(n^2\log n)$. 
This indicates that the presence of community structures can accelerate the anonymization process at least locally, reducing the number of steps needed to obscure the original graph's structure.
Our goal here is to provide a simple example where anonymization occurs faster than mixing; therefore, we leave further explore such sparser graphs in this context for future work.

To validate our theoretical predictions, we performed a series of simulations on both synthetic and real-world networks. The simulations supported our theoretical claims, with regression analysis showing that the slopes of the fitted matchability curves align well with the predicted leading order of the anonymization time. 
We further extended our analysis to real-world networks, including a Facebook friendship network and an email communication network from a European research institution. These experiments provided practical insights into how the proposed models perform on real data, demonstrating the applicability of our theoretical results beyond synthetic settings. In particular, the standard edgelighter walk on the Facebook network and the block edgelighter walk on the email network both showed that anonymization occurred in line with the expected theoretical bounds, underscoring the soundness of our theories.

The findings of this study open new avenues for future research. First, while our results demonstrate that community structures can accelerate anonymization, our current results are only based on edge structure information and it is worth exploring the effect of node similarities in this context. Additionally, extending our models to other types of structured graphs, such as hierarchical stochastic block models \cite{clauset08:_hierar,li2022hierarchical,peixoto14:_hierar_block_struc_high_resol,noroozi2022hierarchy,lyzinski2016community}, block models with overlapping communities \cite{Airoldi2008}, or graphs with more complex structures \cite{Hoff2002,Athreya2018RDPG}, could provide deeper insights into the relationship between graph structure and anonymization time.
Moreover, if we consider collections of temporal networks, dependencies may emerge in the propagation of signals from one graph to the next, extending beyond the scope of Markovian noise. Further analysis of matchability and other statistical inference tasks in the presence of temporal dependence could provide invaluable insights for more realistic real-data modeling via networks. Finally, while the Seeded Graph Matching (SGM) algorithm provided a practical solution for matching noisy graphs to their originals, it would be valuable to explore more advanced uses of seed information---beyond simply constraining the solver to a particular subspace of the function space---to improve the efficiency and accuracy of the matching process; see \cite{mossel2020seeded} for steps in this direction. 

\newpage
\bibliographystyle{plain}
\bibliography{biblio.bib}

\newpage

\section{Appendix}
Herein we collect the proofs of the main results from the paper.

\subsection{Proof of Lemma \ref{lem:covermix}}
\label{app:covermix}

First note that the normalizing constant for the stationary distribution of the standard Edgelighter Markov Chain, $\pi^\circ$, is
\begin{align*}
   \Psi= \sum_{(x,c)}\frac{1}{n}q_2^{\#\text{ of on edges in }c}q_1^{\#\text{ of off edges in }c}=\begin{cases}\frac{q_2^{\binom{n}{2}+1}-\,q_1^{\binom{n}{2}+1}}{q_2-q_1}&\text{ if }q_2\neq q_1;\\
   (\binom{n}{2}+1)q_2^{\binom{n}{2}}&\text{ if }q_2= q_1.
    \end{cases}
\end{align*}
In either case, this normalizing constant is no greater than 1.
We will closely follow the proof of Theorem 19.2 from \cite[Section 19.3]{levin2017markov}; necessary details are provided for completeness.
For the upper bound, let $\vec{0}$ be the configuration with all edges turned off, and $\vec {1}$ the configuration with all edges on.
Let $y$ minimize 
$\mathbb{P}_x(L_t=\cdot\,|\tau^{(u)}_{\text{cover}}\leq t)$ (so that $\mathbb{P}_x(L_t=y|\tau^{(u)}_{\text{cover}}\leq t)\leq 1/n$) and note that
\begin{align*}
\frac{\mathbb{P}_{(x,\vec{0})}(Z_t=(y,\vec{1}))}{\pi^\circ(y,\vec{1})}&=\frac{\mathbb{P}_x(\tau^{(u)}_{\text{cover}}\leq t)q_2^{\binom{n}{2}}\mathbb{P}_x(L_t=y|\tau^{(u)}_{\text{cover}}\leq t)}{\frac{1}{n}q_2^{\binom{n}{2}}/\Psi}\\
&\leq n\Psi\mathbb{P}_x(\tau^{(u)}_{\text{cover}}\leq t)\frac{1}{n}\\
&\leq \mathbb{P}_x(\tau^{(u)}_{\text{cover}}\leq t)
\end{align*}
This yields that 
the separation distance for $(Z_t)$ satisfies
$$s^\circ(t)=\max_{(x,c),(y,c')}1-\frac{\mathbb{P}_{(x,c)}(Z_t=(y,c'))}{\pi^\circ(\,(y,c')\,) }\geq \mathbb{P}_x(\tau^{(u)}_{\text{cover}}> t).$$

Here, we also have that if 
\begin{align*}
d^\circ(t)&=\max_{(x,c)}\|\mathbb{P}^t_{(x,c)}(\cdot)-\pi^\circ\|_{TV}\\ 
\bar d^\circ(t)&=\max_{(x,c),(y,c')}\|\mathbb{P}^t_{(x,c)}(\cdot)-\mathbb{P}^t_{(y,c')}(\cdot)\|_{TV}
\end{align*}
then, as $(Z_t)$ is reversible,
\begin{align*}
d^\circ(t)\leq \bar d^\circ(t)\leq 2d^\circ(t)\\ 
s^\circ(2t)\leq 1-(1-\bar d^\circ(t))^2
\end{align*}
where the last line's inequality follows from Lemma 4.7 in \cite{aldous2002reversible}.
Therefore, as $\bar d^\circ(t)$ is submultiplicative (Lemma 4.11 in \cite{levin2017markov}) and note that by Definition \ref{def:mix}, $d^\circ(t^\circ_m)\leq 1/4$,
$$\bar d^\circ(2t^\circ_m)\leq(\bar d^\circ(t^\circ_m))^2\leq 4  d^\circ(t^\circ_m)^2\leq 1/4$$
so for any $x$,
$$\mathbb{P}_x(\tau^{(u)}_{\text{cover}}> 4t^\circ_m)\leq s^\circ(4t^\circ_m)\leq 1-(1-\bar d^\circ(2t^\circ_m))^2<1/2.$$ 
The Markov property applied successively at each multiple of $4t^\circ_m$ then yields for all integers $k\geq 1$,
$$\mathbb{P}_x(\tau^{(u)}_{\text{cover}}> 4kt^\circ_m)\leq \mathbb{P}_x(\tau^{(u)}_{\text{cover}}> 4t^\circ_m)^k<(1/2)^k.$$
Therefore
\begin{align*}
\mathbb{E}_x(\tau^{(u)}_{\text{cover}})=\sum_{k=0} \mathbb{P}_x(\tau^{(u)}_{\text{cover}}>k)\leq 4t^\circ_m\sum_{k=0} \mathbb{P}_x(\tau^{(u)}_{\text{cover}}>k4t^\circ_m)\leq 8t^\circ_m.
\end{align*}

For the lower bound, again we closely follow the argument in \cite[Section 19.3]{levin2017markov}.  
Let $Z_t=(L_t,C_t)$ be the state of the edgelighter at time $t$.
Set $w=8t^{(u)}_{\text{cover}}+t_m(L_t,1/8)$ (where $\pi_{L}$ is the stationary distribution for $(L_t)$---i.e., the uniform probability on $V$---and $t_m(L_t,1/8)$ is the minimum time such that $\max_x\|\mathbb{P}_x(L_t\in\cdot)-\pi_{L}\|_{TV}\leq 1/8$).
Fix $(v,\vec{0})$ and define 
$$\mu_s=\mathbb{P}_{(v,\vec{0})}(Z_w\in \cdot |\tau^{(u)}_{\text{cover}}=s)$$
We then have that 
$$
\mathbb{P}_{(v,\vec{0})}(Z_w\in \cdot)=\sum_s \mu_s \mathbb{P}_{(v,\vec{0})}(\tau^{(u)}_{\text{cover}}=s),
$$
and so
$$
\|\mathbb{P}_{(v,\vec{0})}(Z_w\in \cdot)-\pi^\circ\|_{TV}\leq \sum_s \mathbb{P}_{(v,\vec{0})}(\tau^{(u)}_{\text{cover}}=s)\|\mu_s-\pi^\circ\|_{TV}
$$
Markov's inequality yields $\mathbb{P}(\tau^{(u)}_{\text{cover}}>8t^{(u)}_{\text{cover}})\leq 1/8$ yielding
$$
\|\mathbb{P}_{(v,\vec{0})}(Z_w\in \cdot)-\pi^\circ\|_{TV}\leq 1/8+\sum_{s\leq8t^{(u)}_{\text{cover}}}  \mathbb{P}_{(v,\vec{0})}(\tau^{(u)}_{\text{cover}}=s)\|\mu_s-\pi^\circ\|_{TV}
$$
If $\tau^{(u)}_{\text{cover}}=s$ for $s\leq w$, then $C_w$ has distribution (for a configuration {$c\in C$}, let $m_c$ and $r_c$ denote the respective number of $1$'s and $-1$'s in $c$)
$$\mathbb{P}(C_w=\vec c)=\nu(c)\propto q_2^{m_c}q_1^{r_c}$$
independent of the location of $L_w$.
Therefore 
\begin{align*}
\mu_s&=\sum_x \mathbb{P}_{(v,\vec{0})}(Z_w\in \cdot |\tau^{(u)}_{\text{cover}}=s,L_s=x)\mathbb{P}_{(v,\vec{0})}(L_s=x|\tau^{(u)}_{\text{cover}}=s)\\ 
&=\sum_x(\mathbb{P}_x(L_{w-s}\in\cdot)\times\nu )\mathbb{P}_{(v,\vec{0})}(L_s=x|\tau^{(u)}_{\text{cover}}=s)
\end{align*}
As $\pi^\circ=\pi_{L}\times \nu$, we have
\begin{align*}
\|\mu_s-\pi^\circ\|_{TV}&\leq 
\sum_x \| \mathbb{P}_x(L_{w-s}\in\cdot)\times\nu-\pi^\circ\|_{TV} \mathbb{P}_{(v,\vec{0})}(L_s=x|\tau^{(u)}_{\text{cover}}=s)\\ 
&\leq \max_x \| \mathbb{P}_x(L_{w-s}\in\cdot)-\pi_L\|_{TV}
\end{align*}
For $s\leq8t^{(u)}_{\text{cover}}$, $w-s>t_m(L_t,1/8)$ and hence $\|\mu_s-\pi^{\circ}\|_{TV}\leq 1/8$.
Therefore,
$$\|\mathbb{P}_{(v,\vec{0})}(Z_w\in \cdot)-\pi^\circ\|_{TV}\leq 1/4,$$
and $w\geq t_m^\circ$.
Submultiplicativity of $\bar d$ yields that
$$t_m(L_t,1/8)\leq 3t_m(L_t)\leq 3\leq 3t^{(u)}_{\text{cover}},$$
where the second inequality follows from Chapter 4 of \cite{aldous2002reversible}, and the fact that $T_j$ is the hitting time of state $j$.
Hence, 
$$
t_m^\circ \leq w\leq 11t^{(u)}_{\text{cover}}
$$
as desired.

\subsection{Derivation of Equations~\eqref{eq:covn2} and \eqref{eq:pijt}}
\label{app:pijt}
    \subsubsection{Proving Eq. \eqref{eq:pijt}}
We seek to show the given bound on $\mathfrak{p}_{t,u,v}$.
We first define the Global Edgelighter Walk analogously to the standard edgelighter, except that the position of the edgelighter (i.e., $L_t$) does not matter.
That is to say, at the $(t+1)$-th step, the global edgelighter randomly chooses any $e\in \binom{V}{2}$, and 
defines $h_{G_{t+1}}(e)$ as in Eqs.\@ \eqref{eq:lamp_on_to_off} and \eqref{eq:lamp_off_to_on}.
Consider the three events, for a given $\{u,v\}\in\binom{V}{2}$:
\begin{align*}
    \mathcal{E}_{1,t}&=\{\text{the edgelighter in the standard walk does not traverse }\{u,v\}\text{ by time }t\}\\
    \mathcal{E}_{2,t}&=\left\{\text{the edgelighter in the standard walk does not traverse }\{u,v\}\text{ in }\big(\{L_{2i},L_{2i+1}\}\big)_{i=0}^{\lfloor \frac{t-1}{2}\rfloor}\right\}\\
    \mathcal{E}_{3,t}&=\{\text{the edgelighter in the global walk does not traverse }\{u,v\}\text{ by }t\}
\end{align*}
We then have
$$\left(1-\frac{1}{\binom{n}{2}}\right)^{t}=\PP(\mathcal{E}_{3,t})\leq \PP(\mathcal{E}_{1,t})\leq \PP(\mathcal{E}_{2,t})=\left(1-\frac{1}{\binom{n}{2}+n}\right)^{\lfloor \frac{t+1}{2}\rfloor}$$
Using the fact that for $n>3$, 
$$
\left(1-\frac{1}{\binom{n}{2}+n}\right)\leq e^{-\frac{1}{\binom{n}{2}+n}}
~\text{ and }~
\left(1-\frac{1}{\binom{n} {2}}\right)\geq e^{-\frac{1}{\binom{n}{2}-1}},
$$
we have that 
\begin{align*}
    \mathfrak{p}_{t,u,v}&\leq \text{exp}\left\{ -\frac{\lfloor \frac{t+1}{2}\rfloor}{\binom{n}{2}+n}\right\}
    ~\text{ and }~
\mathfrak{p}_{t,u,v}\geq\text{exp}\left\{ -\frac{t}{\binom{n}{2}-1}\right\},
\end{align*}
which establishes Equation~\eqref{eq:pijt}.

\subsubsection{Proving Eq. \eqref{eq:covn2}}
    Letting $\tau_{\text{cover}}^{(g)}$ be the cover time of the global edgelighter, and $\tau_{\text{cover}}$ the cover time defined in Eq. \ref{eq:cov}, we have that $\tau_{\text{cover}}^{(g)}\leq_{st.}\tau_{\text{cover}}^{(u)}\leq_{st.}2\tau_{\text{cover}}$ (where $\leq_{st.}$ denotes stochastic ordering). 
    We therefore have that for $n$ sufficiently large, via simple coupon collector asymptotics,
    \begin{align*}
\frac{1}{2}n^2\log n\leq \mathbb{E}\tau_{\text{cover}}^{(g)}\leq \max_x\mathbb{E}_x\tau_{\text{cover}}^{(u)}\leq \max_x\mathbb{E}_x\tau_{\text{cover}}\leq \frac{5}{2} n^2\log n,
    \end{align*}
    which yields Equation~\eqref{eq:covn2}.

\subsection{Proof of Theorem \ref{thm:global}}
\label{app:thm1}

Here, we will follow standard concentration inequality theory to establish the matchability result.
Let $B=A_t$ and $A=A_0$, and for a fixed permutation matrix $P$, with associated permutation $\sigma_P$, consider the quantity
    \begin{align*}
    \snp &:= \tr{APBP^T} -\tr{AB} = \tr{ A[PBP^T-B] } \\
    &= 2\sum_{\substack{\{i,j\}\in\binom{V}{2}\text{ s.t. }\\\{i,j\}\neq\{\sigma_P(i),\sigma_P(j)\}} } A_{i,j} \left(\bsij-\bij\right).
    \end{align*}
Straightforward computations then yield    
    \begin{align*}
        \E \snp &= 2\sum_{\substack{\{i,j\}\in\binom{V}{2}\text{ s.t. }\\\{i,j\}\neq\{\sigma_P(i),\sigma_P(j)\}} } \E\left\{\left(\bsij - \bij\right)\mid \aij=1\right\}\PP(\aij=1)\\
        &= 2p\sum_{\substack{\{i,j\}\in\binom{V}{2}\text{ s.t. }\\\{i,j\}\neq\{\sigma_P(i),\sigma_P(j)\}} } \left\{\E\left(\bsij\mid \aij=1\right) - \E\left(\bij\mid \aij=1\right)\right\}
    \end{align*}
Noting that $\mathfrak{p}_{t,i,j}$ does not depend on $i,j$ in this case, we drop these indices and write $\mathfrak{p}_{t}$ moving forward.
    We now consider computing each conditional expectation via considering cases for index $i$ and $j$ after $t$ steps.
    \begin{enumerate}
    \item We first consider $\E\left\{\bsij \mid \aij=1\right\}$.  
    If $\{\sigma_P(i),\sigma_P(j)\}$ was traversed by the edgelighter by time $t$, then the edge has probability $p$ regardless of the state of $A_{\{\sigma_P(i),\sigma_P(j)\}}$ and $\aij$.  This conditional expectation is then equal to
    $$\E\left\{\bsij \mid \aij=1\right\}=p(1-\ppt)+p\ppt=p$$
    
      \item We next consider $\E\left\{\bij \mid \aij=1\right\}$.  
    If $\{i,j\}$ was traversed by the edgelighter by time $t$, then the edge has probability $p$ regardless of the state of $\aij$; else $\bij=\aij$.  This conditional expectation is then equal to
    $$\E\left\{\bij \mid \aij=1\right\}=p(1-\ppt)+\ppt$$
\end{enumerate}

Putting the above together, we get that 
\begin{align*}
    \E(\snp)=2\underbrace{\left|\left\{ \{i,j\}\in\binom{V}{2}\text{ s.t. } \{i,j\}\neq \{\sigma_P(i),\sigma_P(j)\}\right\} \right|}_{=:\eta_P}\ppt p(p-1)
\end{align*}
If $P$ shuffles $k$ vertices, then $\eta_P\in [(n-2)k/2 ,nk]$ and this expectation can be bounded as
$$
(n-2)k\ppt p(p-1)\leq \E(\snp)\leq 2nk\ppt p(p-1) .
$$
Moreover, $\snp$ can be realized as a function of the following random variables:
\begin{itemize}
\item[i.] The  at most $2nk$ edges, $\aij$, involved in the summation. Changing any of these could change the summation $\snp$ by at most 4.
\item[ii.] The  at most $2(t+1)$ random variables dictating the steps and edge flips of the edgelighter walk.
Changing a single move of the edgelighter would change at most $2$ of the $\bij$'s and hence could change the summand by at most 8.  Changing a flip of the edgelighter could change at most one $\bij$ and hence would change the summand by at most 2.
\end{itemize}
McDiarmid's inequality then gives us that for $n$ sufficiently large (recalling that we assume there exists a constant $c>0$ such that $t\leq c n^2\log n$)
\begin{align*}
\PP(\snp\geq 0)&\leq \PP(|\snp-\E(\snp)|\geq |\E\snp|)\\
&\leq 2\text{exp}\left\{ -\frac{(n-2)^2k^2 p^2(p-1)^2\text{exp}\left\{ -\frac{4t}{n(n-1)-2}\right\} }{16(nk+t+1)} \right\}\\
&\leq 2\text{exp}\left\{ -\frac{(n-2)^2k^2 p^2(p-1)^2n^{-5c} }{16(nk+c n^2\log n+1)} \right\} .
\end{align*}
Let $\Pi_{n,k}$ be the set of permutations shuffling exactly $k$ vertex labels, and define the event
\begin{equation*}
\beta_{k,t}=\{\exists P\in\Pi_{n,k}\text{ s.t. }\snp\geq 0\}.
\end{equation*}
A union bound over the at most $n^k$ permutation matrices in $\Pi_{n,k}$ yields 
\begin{align*}
\PP(\beta_{k,t})&\leq  2\exp\left\{ -\frac{(n-2)^2k^2 p^2(p-1)^2n^{-5c} }{16(nk+c n^2\log n+1)}+k\log n \right\} .
\end{align*}
Recalling that by assumption there exists a constant $\alpha>5c$ such that $k\geq n^\alpha$, we have 
\begin{align*}
\frac{1}{k\log n} ~\frac{(n-2)^2k^2 p^2(p-1)^2n^{-5c} }{16(nk+c n^2\log n+1)}
=\Omega\left(\frac{p^2(p-1)^2n^{\alpha-5c}}{\log^2 n}  \right)=\omega(1) .
\end{align*}
Therefore, taking a union over such $k$ yields
\begin{align*}
\PP\left(\bigcup_{k\geq n^\alpha}\beta_{k,t}\right)&\leq  2\text{exp}\left\{ -\Omega\left(\frac{n^{2\alpha-5c} }{\log n}\right) +\log n\right\}\leq  2\text{exp}\left\{ -\Omega\left(\frac{n^{2\alpha-5c} }{\log n}\right)\right\},
\end{align*} 
completing the proof.
\vspace{3mm}

\subsection{Proof of Theorem \ref{thm:globalbad}} \label{sec:thm2}

Abusing notation slightly, we will use $\Pi_n$ to denote both the set of permutations of $\{1,2,\cdots,n\}$ (denoted by lower case Greek letters) and the set of $n\times n$ permutation matrices (denoted by capital Roman letters), with the intended meaning made clear from context. Without loss of generality, we assume throughout this proof that $n^\beta$ is an integer, as otherwise we could simply replace $n^\beta$ with $\lfloor n^\beta\rfloor$ throughout the argument that follows.
We say that permutations $\sigma,\tau \in \Pi_n$ disagree at location $\ell\in\{1,2,\cdots,n\}$ if $\sigma(\ell)\neq \tau(\ell)$. Consider $\mathfrak{n}_\beta=\lfloor n/(2n^\beta +1)\rfloor$ disjoint sets in $\{1,2,\cdots,n\}$, each with size $2n^\beta +1$ (except for the last set, which may contain a few more elements). We denote these sets via $\{S_i\}_{i=1}^{\mathfrak{n}_\beta}$.
For each $i\in\{1,2,\cdots,\mathfrak{n}_\beta\}$, let $\sigma^{(i)}$ be a permutation of $\{1,2,\cdots,n\}$ that fixes all elements of $\cup_{j : j\neq i}S_j$ and is (in cycle decomposition) a cyclic derangement (i.e., a cycle of length $2n^\beta+1$ of $S_i$, except for the potentially longer cycle of $S_{\mathfrak{n}_\beta}$). 

Note that if $(\sigma^{(i)})^k$ denotes the permutation $$(\sigma^{(i)})^k=\underbrace{\sigma^{(i)}\circ\cdots\circ \sigma^{(i)}}_{k\text{ compositions }},$$
then for all $k_1,k_2\in\{1,2,\cdots,2n^\beta\}$ with $k_1\neq k_2$, we have that $(\sigma^{(i)})^{k_1}$ and $(\sigma^{(i)})^{k_2}$ disagree in at least $2n^\beta +1$ locations (and exactly $2n^\beta+1$ locations for all $i\neq \mathfrak{n}_\beta$).
Consider all permutations of the form
\begin{align}
\label{eq:permperm}
\sigma=(\sigma^{(1)})^{k_1}\circ(\sigma^{(2)})^{k_2}\circ \cdots \circ (\sigma^{(\mathfrak{n}_\beta)})^{k_{\mathfrak{n}_\beta}}
\end{align}
for $(k_1,k_2,\cdots,k_{\mathfrak{n}_\beta})\in\{1,2,\cdots,2n^\beta\}^{k_{\mathfrak{n}_\beta}}$.
Each pair of such $\sigma$'s with distinct $(k_1,k_2,\cdots,k_{\mathfrak{n}_\beta})$ sequences (i.e., there is at least one $j$ such that their $j$-th entries differ) disagree in at least $2n^\beta+1$ locations.
We denote by $\widetilde{\Pi}^{\beta}$ the set of all permutations of the form in Eq.~\eqref{eq:permperm} with $(k_1,k_2,\cdots,k_{\mathfrak{n}_\beta})\in{\{1,2,\cdots,2n^\beta\}^{k_{\mathfrak{n}_\beta}}}$.
Note that the size of $\widetilde{\Pi}^{\beta}$ is 
\begin{equation} \label{eq:Pitilde:count}
\left |\widetilde{\Pi}^{\beta}\right|=(2n^\beta)^{\mathfrak{n}_\beta}\geq (2n^\beta)^{n/(2n^\beta+1)-1}.
\end{equation}

For a set of permutations $\mathcal{Q}\subset\Pi_n$, define the action of a permutation $P\circ \mathcal{Q}$ to be the set
\begin{equation*}
P\circ \mathcal{Q}=\{Q\in\Pi_n\text{ s.t. }P^{-1}Q\in \mathcal{Q}\} .
\end{equation*}
For a permutation $P \in \Pi_n$, define the event
\begin{equation*}
\mathcal{E}_{\beta,P}=\left\{\argmax_{Q\in\Pi_n} \tr{A_0QA_tQ^T}\subset P\circ\left(\bigcup_{k=0}^{n^\beta}\Pi_{n,k}\right)  \right\}.
\end{equation*}
Note that the sets $\{\mathcal{E}_{\beta,P}\}_{P\in\widetilde{\Pi}^{\beta}\cup\{I_n\}}$ are disjoint, as the sets $P\circ\left(\bigcup_{k=0}^{n^\beta}\Pi_{n,k}\right)$ are disjoint for $P\in\widetilde{\Pi}^{\beta}\cup\{I_n\}$.
To see that this is the case, suppose by way of contradiction that for distinct $P_1,P_2\in \widetilde{\Pi}^{\beta}\cup\{I_n\}$, there is a $Q\in P_1\circ\left(\bigcup_{k=0}^{n^\beta}\Pi_{n,k}\right)\cap P_2\circ\left(\bigcup_{k=0}^{n^\beta}\Pi_{n,k}\right)$.
Then there exist $Q_1,Q_2\in \bigcup_{k=0}^{n^\beta}\Pi_{n,k}$ such that 
$P_1^{-1}Q=Q_1$ and $P_2^{-1}Q=Q_2$, and hence $P_1=P_2Q_2Q_1^{-1}$.
This would imply $P_1$ and $P_2$ could disagree in at most $2n^\beta$ locations, but this contradicts our assumption that $P_1,P_2\in \widetilde{\Pi}^{\beta}\cup\{I_n\}$ are distinct, since distinct elements of $\widetilde{\Pi}^{\beta}\cup\{I_n\}$ must differ in at least $2n^\beta +1$ locations by construction.

Recall at the beginning of Section \ref{sec:Main}, we mentioned that for $t\geq \tau^{(u)}_{\text{cover}}$, $A_t\sim \ER(n,p)$ is independent of $A_0$ and hence $\tr{A_0PA_tP^T}$ is equal in distribution to $\tr{A_0A_t}$ for any permutation $P\in\Pi_n$.
It follows that $\mathbb{P}(\mathcal{E}_{\beta,P}\,|\,\,t > \tau^{(u)}_{\text{cover}})$ is the same for all $P\in \Pi_n$.
Therefore, in light of Equation~\eqref{eq:Pitilde:count},
\begin{equation} \label{eq:bound:PcalEmidCover}
  \mathbb{P}(\mathcal{E}_{\beta,I_n}\,|\,\,\,t > \tau^{(u)}_{\text{cover}}) \leq (2n^\beta)^{1-n/(2n^\beta+1)}.
\end{equation}
Theorem~\ref{thm:globalbad} is then proven by considering
\begin{align*}
\PP(\mathcal{E}_{\beta,I_n}) = &\PP(\mathcal{E}_{\beta,I_n}|\tau^{(u)}_{\text{cover}}\geq t)\PP(\tau^{(u)}_{\text{cover}}\geq t) + \PP(\mathcal{E}_{\beta,I_n}|t > \tau^{(u)}_{\text{cover}})\PP(t > \tau^{(u)}_{\text{cover}})\\
\leq& \PP(\tau^{(u)}_{\text{cover}}\geq t)+\PP(\mathcal{E}_{\beta,I_n}|\,\,t > \tau^{(u)}_{\text{cover}})
\end{align*}
and applying Equation~\eqref{eq:bound:PcalEmidCover} along with our bound in Equation~\eqref{eq:couponcollector}, using our assumption that $t \ge 4 n^2 \log n $.

\subsection{Proof of Theorem \ref{thm:SBM}}
\label{sec:pf3}

For each $i\in[K_n]$, let $n_i$ denote the size of the $i$-th community.
We will make use of the following Chernoff-style bound, adapted from Theorem 3.2 of \cite{chung2006concentration}):
Let $X_i\sim \Bern(p_i)$ be independent random variables with $X=\sum_i X_i$ and $\mathbb{E}(X)=\sum_i p_i$.
Then for any $t>0$, 
\begin{align*}
    \mathbb{P}\left(|X-\mathbb{E}(X)|>t\right) \leq 2 \exp\left\{ -\frac{t^2}{2\mathbb{E}(X)+2t/3}\right\} .
\end{align*}
Applying this to $m_{ij}\sim \Bin(n_in_j,\Lambda_{i,j} )$ defined to be the number of edges between community $i$ and community $j$, with $m_{ii}\sim \Bin(\binom{n_i}{2},\Lambda_{i,i} )$ defined analogously, we have that with probability at least $1-K_n^2 n^{-3}$ (as $n_in_j\Lambda_{i,j}\geq 5\log n$ by assumption),
\begin{align}
\label{eq:across}
    \left| m_{ij}-n_in_j\Lambda_{i,j}  \right|&\leq\sqrt{5n_in_j\Lambda_{i,j}\log n} ~\text{ and }\\
    \label{eq:within}
    \left| m_{ii}-\binom{n_i}{2}\Lambda_{i,i}  \right|&\leq\sqrt{5\binom{n_i}{2}\Lambda_{i,i}\log n}.
\end{align}
For the remainder of the proof, we will condition on the events in Equations~\eqref{eq:across} and~\eqref{eq:within}, calling this event $\mathcal{A}_n$, so that between community $1$ and community $j$, there are $m_{1,j}\in n_1n_j\Lambda_{1,j} \pm\sqrt{5n_1n_j\Lambda_{1,j}\log n}$ edges.
Recall that $\Lambda_{1,j}\leq b_3 n^{a_2}/n_1n_j$, so that if the above events hold, there exists a  constant $b_3'>0$ with $m^*=\max_j m_{1,j}\leq b_3'n^{a_2}$.

Consider the edgelighter walk on the random graph on event $\mathcal{A}_n$.
If the walker is in community $1$, then with probability $1/2$ the edgelighter will move to a uniformly random community other than community $1$ on the next step.
Denoting the number of steps before leaving community $1$ by $T_1$, note that 
\begin{equation*}
\mathbb{P}\left(T_1>K_n m^*\right)=\left( \frac{1}{2} \right)^{K_n m^*}\leq \exp\{ -m^* \},
\end{equation*}
and hence $T_1$ is, with high probability, 
of lower order than $N=n^{1+a_2}\log^{3/2}n$ as defined in the theorem. 
As such, $T_1$ has little impact on the overall mixing time, and we can assume without loss of generality that the edgelighter starts outside of community 1.

Once the edgelighter first leaves community 1 (take this to be time 0 if the edgelighter starts outside community 1), we define a renewal process as follows:
We say that a renewal occurs if the edgelighter enters community $1$, immediately moves to a different vertex inside community $1$, and then eventually exits community $1$.
The renewal is defined to occur when all the conditions are met, i.e., when the edgelighter eventually exits community $1$. 
Below, we write $X\sim \Geo(p)$ if $\mathbb{P}(X=k)=(1-p)^{k} p$ for $k=0,1,2,\dots$. 
The inter-arrival time for this renewal process then has distribution equal to
\begin{equation*}
\xi_n\stackrel{\text{dist.}}{=}\sum_{i=1}^{Z_1}(V_i+W_i+2)+(Z_2+3),
\end{equation*}
where
\begin{itemize}
    \item[i.] $\{V_i\}$ are i.i.d. $\Geo\left(\frac{1}{2(K_n-1)}\right)$ random variables corresponding to the amount of time spent outside community $1$ between the $(i-1)$-th and $i$-th visits to community $1$;
    \item[ii.] $\{W_i\}$ are i.i.d.~and count the number of moves inside community 1 during a non-renewal visit (i.e., where the edgelighter does not move to a new vertex in community 1 in its first move) such that we know $W_i=0$, if the edgelighter immediately leaves community 1 (which happens with probability $1/2$); and $W_i=1+\widetilde W_i$, if the edgelighter stays at the current node for the immediate next step (which happens with probability $(1/2)\cdot (1/n)$), where $\widetilde W_i\sim \Geo(1/2)$ (note the summand $+2$ is included to count the steps entering and exiting community 1);
    \item[iv.] $Z_1\sim\Geo(1/2(1+1/n))$, independent of $\{V_i,W_i\}$, represents the number of non-renewal visit to community $1$; i.e., where the first move is to exit community 1 or is a ``hold'' (i.e., stays for the immediate next step upon entering community 1);
    \item[v.] $Z_2\sim\Geo(1/2)$, independent of $\{V_i,W_i\}$ and $Z_1$, represents the number of moves inside community $1$ (not including the initial move to a new vertex) the edgelighter takes in the final visit to community $1$ before renewal (note the summand $+3$ term counts the step entering community 1, the immediate move to a new vertex in community 1 and the step exiting community 1).
\end{itemize}
Here ,
\begin{equation*}
\mathbb{E}(\xi_n)=\frac{2n}{n+1}\left(2(K_n-1)+\frac{3}{2n} +2\right)+5\leq 12 K_n .
\end{equation*}
Let $R(t)$ denote the number of renewals by time $t$ and denote the $m$-th renewal time via $S_m$ for $m = 0,1,2,\dots$.
Let 
$N = n^{1+a_2}\log^{3/2}(n)$ and 
$f(n) = n^{(3a_2+1)/2}\log^{3/2}(n)$, 
where we recall that $0<a_2<1$ is the constant picked by the theorem to regulate the edge probability between community 1 and other communities.
By Markov's inequality, there exists a constant $C>0$ such that
\begin{align*}
\mathbb{P}\left(R(N)< f(n) \right)&=\mathbb{P}\left(S_{f(n)} > N\right)
\leq \frac{\mathbb{E}(\xi_n)n^{(3a_2+1)/2}\log^{3/2}(n)}{ n^{1+a_2}\log^{3/2}(n)} 
\leq C n^{(a_2-1)/2}.
\end{align*}

Between renewals, the edgelighter must enter community 1 at least once.
This entry is from a uniformly random community, {as the random walk at the community level is the simple symmetric random walk on the $K$ communities.  Hence, once the edgelighter leaves community 1 (which must be done between renewals), the edgelighter would re-enter community 1 from a uniformly random other community.
Moreover, a randomly selected edge and non-edge between the entering community and community 1 are switched.}
When moving from community $i$ to $1$, the randomization of the edges/non-edges between the communities is equivalent to a random walk on the Johnson graph $J(n_1n_i,m_{i1})$ \cite{brouwer2011distance}.
Properties of the spectrum of this random walk are well-known.
In particular, the relaxation time, here defined as $1/(1-\lambda_2)$ \citep{levin2017markov}, where $\lambda_2$ is the second largest eigenvalue of the random walk adjacency matrix, is equal to 
$\trel^{J,i}=m_{i1}(n_in_1-m_{i1})/n_in_1$.
Moreover, the mixing time for this random walk can be bounded as follows:
since the stationary distribution for the random walk on the Johnson graph is uniform and 
\begin{equation*}
\binom{n_in_1}{m_{1i}}\leq \frac{( e n_1n_i)^{m_{1i}}}{ m_{1i}^{m_{1i}}},
\end{equation*}
we have that for $n$ sufficiently large (where $C$ is a constant that can change line-to-line),
\begin{align*}
t_{\text{mix}}^{J,i}(\epsilon)
&\leq \trel^{J,i} \log\left(\frac{( e n_1n_i)^{m_{1i}}}{\epsilon m_{1i}^{m_{1i}}}\right)
\leq m_{1i}\log\left(\frac{( e n_1n_i)^{m_{1i}}}{ \epsilon m_{1i}^{m_{1i}}}\right) \\
&\leq C m_{1i}\log\left(\frac{( n_1n_i)^{m_{1i}}}{ \epsilon }\right)
\leq C m^*\log\left(\frac{( n_1n_i)^{m^*}}{ \epsilon }\right) .
\end{align*}
Let
$t_{\text{mix}}^{J*}(\epsilon)=\max_i t_{\text{mix}}^{J,i}(\epsilon)$.
Setting 
\begin{equation*}
\epsilon=\left(n_1 n^*\right)^{-b_3'n^{a_2}},
\end{equation*}
we have that (again letting $C$ be a constant that can change line-to-line)
\begin{align*}
t_{\text{mix}}^{J*}(\epsilon)=t_{\text{mix}}^{J*}\left(\left(n_1 n^*\right)^{-b_3'n^{a_2}}\right)
&\leq C n^{2a_2}\log\left(n_1n^*\right)\leq C n^{2a_2}\log\left(n\right) .
\end{align*}
Once there have been $t_{\text{mix}}^{J*}(\epsilon)$ transitions from each other community (i.e., all of communities $2, 3, \ldots, K_n$) to community 1, then the signal between community 1 and all other communities has been uniformly mixed (at level $\epsilon$).
Standard coupon collector asymptotics then give us that, for an appropriate constant $C$ that can change from line to line, after
\begin{equation*}
C n^{(1-a_2)/2}n^{2a_2}\log\left(n\right)K_n\log K_n\leq C n^{(3a_2+1)/2}\log\left(n\right)
\end{equation*}
renewals, we have visited each inter-community network (here the inter-community networks are the $K-1$ bipartite graphs consisting of edges/non-edges between community 1 and community $j\neq 1$) connecting community $1$ at least $t_{mix}^{J*}(\epsilon)$ times with probability at least
\begin{equation*}
    1-C n^{2a_2}\log\left(n\right)b_4^{-n^{(1-a_2)/2}}.
\end{equation*}
Recall that $\epsilon=\left(n_1 n^*\right)^{b_3'n^{a_2}}$, and for each integer $\eta\in\mathbb{Z}>0$, define the events 
\begin{align*}
\mathcal{A}_\epsilon(\eta)=\big\{& \text{For each }i\neq 1,\text{ there have been at least }t_{\text{mix}}^{J*}(\epsilon)\text{ transitions }\\&~~~\text{from community }i\text{ to community } 1 \text{ by time }\eta\big\}
\end{align*}
and
\begin{align*}
\mathcal{B}(\eta)=\big\{&\text{The lamplighter walk has covered the edge-set within community 1 by time }\eta\big\} .
\end{align*}
At a time $\eta$, we say that the edge configurations within community $1$ and between community 1 and all other communities are uniformly random if the following holds:
{\begin{itemize}
    \item[i] Letting $B_1$ be the set of vertices in community 1, we have that $A_\eta[B_1]\sim \ER(n_1,q)$ independent of $A_0[B_1]$; here $A_0[B_1]$ (resp., $A_\eta[B_1]$ is the induced subgraph of $A_0$ (resp., $A_\eta$) on vertices $B_1$);
    \item[ii.] For $i\neq 1$, let $B_i$ be the vertices in community $i$.  Letting $A_0[B_1; B_i]$ (resp., $A_\eta[B_1; B_i]$ be the induced bipartite subgraph of $A_0$ (resp., $A_\eta$) between vertex sets $B_1$ and $B_i$, we have that $A_\eta[B_1; B_i]\sim \ER(n_1,n_i,m_{1i})$ independent of $A_0[B_1;B_i]$; here $\ER(n_1,n_i,m_{1i})$ is the uniform distribution on bipartite graphs between $n_1$ and $n_i$ vertices with $m_{1i}$ edges.  Moreover, the graphs $A_\eta[B_1]$ $\{A_\eta[B_1; B_i]\}_{i\neq 1}^K$ are independent.
\end{itemize}}
\noindent If this is the case, then writing $\buildrel \mathcal{D} \over =$ for equality in distribution, we have that
\begin{equation*}
\tr{A_0(Q_1\oplus I_{n-n_1})A_\eta(Q_1\oplus I_{n-n_1})^T}\buildrel \mathcal{D} \over = \tr{A_0 (Q_2\oplus I_{n-n_1})A_\eta(Q_2\oplus I_{n-n_1})^T} \,\; \forall Q_1, Q_2 \in \Pi_{n_1}
\end{equation*}
Letting the uniform distribution of the edges within and and across community 1 be denoted $\psi=\otimes_{i=1}^{K_n} \psi_i$---$\psi_1$ being the distribution within community 1 and $\psi_i$ the distribution across community 1 and community $i$---it follows that, as in the proof of Theorem 2,
\begin{equation*}
\mathbb{P}_{\psi}(\mathcal{E}_{\eta,\beta})\leq (2n_1^\beta)^{{1}-n_1/(2n_1^\beta+1)}.
\end{equation*} 
Let 
$\mathcal{C}$ is the set of configurations of edges within community 1 and between community 1 and all other communities and let
$$R^{(\eta)}=(R^{(\eta)}_1,R^{(\eta)}_2,\cdots,R^{(\eta)}_{K_n})$$
be the random variable taking values over such configurations, where $R^{(\eta)}_1$ is the configuration within community 1, and $R^{(\eta)}_i$ is the configuration between community 1 and community $i$ for $i>1$.
We then have that
\begin{align*}
\mathbb{P}(\mathcal{E}_{\eta,\beta}|\mathcal{A}_\epsilon(\eta),\mathcal{B}(\eta))=\sum_{\vec c~\in~\mathcal{C}}\mathbb{P}\left[\mathcal{E}_{\eta,\beta}|R=c,\mathcal{A}_\epsilon(\eta),\mathcal{B}(\eta)\right]\mathbb{P}\left[R=c|\mathcal{A}_\epsilon(\eta),\mathcal{B}(\eta)\right] .
\end{align*}
Now, conditioning on $\mathcal{A}_\epsilon(\eta),\mathcal{B}(\eta)$, we have that $\eta$ is greater than the sum
$$
\eta\geq \tau_{cov,1}+\sum_{i=2}^{K_n}t_{mix}^{J,i}(\epsilon),
$$
where $\tau_{cov,1}$ is the cover time of the set of unordered pairs within community 1.
{If we condition on the number of steps taken by the edgelighter within community 1 (call this $\eta_1$) and between community 1 and community $i$ (call this $\eta_{1,i}$) for all $i\neq 1$ satisfying
\begin{align}
\label{eq:enoughsteps}
\eta_1\geq \tau_{cov,1};\quad \text{ and }\quad
\eta_{1,i}\geq t_{mix}^{J,i}(\epsilon)\text{ for each $i\neq 1$};
\end{align}
then the coordinates of $R$ are independent.}
Moreover, writing $P_R$ for the distribution of $R$ we have
\begin{align*}
\|P_R-\psi\|_{TV}=\|\otimes_{i=1}^{K-n} P_{R_i}-\otimes_{i=1}^{K-n}\psi_i\|_{TV}\leq \sum_{i=1}^{K_n}\|P_{R_i}-\psi_i\|_{TV}\leq K_n \epsilon .
\end{align*}
{As the total variation distance is non-increasing here, and $\|P_{R_1}-\psi_1\|_{TV}$ is the distance from stationarity for the random walk within community 1, 
and for $i\neq 1$, $\|P_{R_i}-\psi_i\|_{TV}$ is the distance from stationarity for the random walks (Johnson graphs) between community 1 and community $i$, we have that
\begin{align*}
\|P_R-\psi\|_{TV}=\leq \sum_{i=1}^{K_n}\|P_{R_i}-\psi_i\|_{TV}\leq K_n \epsilon .
\end{align*}
As this bound is uniform over all conditionings satisfying Eq. \ref{eq:enoughsteps} (i.e., on the number of steps taken within and across community 1), we have that
$$
|\mathbb{P}(R=c|\mathcal{A}_\epsilon(\eta),\mathcal{B}(\eta))-\psi(c)|\leq K_n \epsilon,
$$
so that 
\begin{align*}
\mathbb{P}(\mathcal{E}_{\eta,\beta}|\mathcal{A}_\epsilon(\eta),\mathcal{B}(\eta))&=\sum_{\vec c\in\mathcal{C}}\mathbb{P}(\mathcal{E}_{\eta,\beta}|R=c,\mathcal{A}_\epsilon(\eta),\mathcal{B}(\eta))\mathbb{P}(R=c|\mathcal{A}_\epsilon(\eta),\mathcal{B}(\eta))\\
&\leq \sum_{\vec c\in\mathcal{C}}\mathbb{P}(\mathcal{E}_{\eta,\beta}|R=c)\psi(c)+K_n \epsilon\\
&\leq \mathbb{P}_{\psi}(\mathcal{E}_{\eta,\beta})+K_n \epsilon\\
&\leq (2n_1^\beta)^{1-n_1/(2n_1^\beta+1)}+K_n \epsilon .
\end{align*}
Now, letting $\eta=N=n^{1+a_2}\log^{3/2}(n)$, and $\epsilon=(n_1n^*)^{-b_3'n^{a_2}}$ we have that
\begin{align*}
  \mathbb{P}(\mathcal{B}(N))&\leq   
  \mathbb{P}(\mathcal{B}(N)|R(N)\geq f(n))+\mathbb{P}(R(N)< f(n))\\
  &\leq   
  \mathbb{P}(\tau_{cov,1}>f(n))+Cn^{a_2-1}/2\\
  &\leq e^{-\frac{n^{(3a_1+1)}\log^{3/2}n}{n^{2a_1}\log(n^{2a_1})}+1}+Cn^{a_2-1}/2\\
  &=e^{-\Omega(n^{a_1}\log ^{1/2}(n))}+O(n^{a_2-1}/2)),
\end{align*}
and 
\begin{align*}
  \mathbb{P}(\mathcal{A}_\epsilon(N))&\leq   
  \mathbb{P}(\mathcal{A}_\epsilon(N)|R(N)\geq f(n))+\mathbb{P}(R(N)< f(n))\\
  &O( n^{2a_2}\log\left(n\right)b_4^{-n^{(1-a_2)/2}})+O(n^{a_2-1}/2)) .
\end{align*}
Combining the above, we have that
\begin{align*}
  \mathbb{P}(\mathcal{E}_{N,\beta})&\leq
  \mathbb{P}(\mathcal{E}_{N,\beta}|\mathcal{A}_\epsilon(N),\mathcal{B}(N))  +\mathbb{P}(\mathcal{A}_\epsilon(N))+\mathbb{P}(\mathcal{B}(N))\\
  &\leq (2n_1^\beta)^{{1}-n_1/(2n_1^\beta+1)}+K_n(n_1n^*)^{-b_3'n^{a_2}}+O( n^{2a_2}\log\left(n\right)b_4^{-n^{(1-a_2)/2}})\\
  &\quad \quad+O(n^{a_2-1}/2))+e^{-\Omega(n^{a_1}\log ^{1/2}(n))},
\end{align*}
as desired.

\subsection{Proof of the stationary distribution for the block edgelighter}
\label{app:stat2}

Here we show that the form of the stationary distribution is
\begin{equation*}
\pi^\bullet(B_i,u,\vec c)\propto\frac{1}{K} \frac{1}{|B_i|}\prod_{j}q_{j;2}^{m_{j,\vec{c}}}q_{j;1}^{r_{j,\vec{c}}} \frac{1}{\prod_{i,j}\binom{n_in_j}{m_{ij}}}
\end{equation*}
by showing that this distribution satisfies detailed balance.
We consider the following cases:
\begin{itemize}
\item[Case 1:] $u\neq v\in B_i$ for some $B_i$
\begin{itemize} 
\item[i.] $\vec{c}_v=\vec{c}_u$ in which case $\pi^\bullet(B_i,u,\vec{c}_u)=\pi^\bullet(B_i,v,\vec{c}_v)$ and we also have 
\begin{align*}
    \mathbb{P}_{(B_i,u,\vec{c}_u)}&(W_1=(B_i,v,\vec{c}_v)=\mathbb{P}_{(B_i,v,\vec{c}_v)}(W_1=(B_i,u,\vec{c}_u)\\
    &=\begin{cases} \frac{1}{2}\frac{1}{|B_i|}(1-q_{i,1})\text{ if }h_{G_t}(\{u,v\})=1\\
\frac{1}{2}\frac{1}{|B_i|}(1-q_{i,2})\text{ if }h_{G_t}(\{u,v\})=\trz,
\end{cases}
\end{align*}
so that detailed balance holds.
\item[ii.] $\vec{c}_1\neq \vec{c}_2$; without loss of generality let the $\{u,v\}$ entry of $\vec c_u$ be $\trz$ and of $\vec c_v$ be $1$ (the opposite case is analogous and left to the reader). 
In this case,
$$\frac{\pi^\bullet(B_i,u,\vec{c}_u)}{\pi^\bullet(B_i,v,\vec{c}_v)}=\frac{q_{i;1}}{q_{i;2}} .
$$
Now
\begin{align*}
    &\mathbb{P}_{(B_i,u,\vec{c}_u)}(W_1=(B_i,v,\vec{c}_v))= \frac{1}{2}\frac{1}{|B_i|}q_{i,2};\\
    &\mathbb{P}_{(B_i,v,\vec{c}_v)}(W_1=(B_i,u,\vec{c}_u))=\frac{1}{2}\frac{1}{|B_i|}q_{i,1}
\end{align*}
and again detailed balance follows immediately.
\end{itemize}
\item[Case 2:] $u\in B_i$,  $v\in B_j$ for some $B_j\neq B_i$.
Wlog, let $\vec c_u$ and $\vec c_v$ differ only at $e,e'$ between blocks $B_i$ and $B_j$, where $h_{\vec c_u}(e)=1$, $h_{\vec c_u}(e)=\trz$, and $h_{\vec c_v}(e)=\trz$, $h_{\vec c_v}(e)=1$. 
To go from $\vec c_u$ to $\vec c_v$, the edge at $e$ is selected to be removed and the non-edge at $e'$ to be added, and vice versa. Thus,
\begin{align*}
    &\mathbb{P}_{(B_i,u,\vec{c}_u)}(W_1=(B_j,v,\vec{c}_v))= \frac{1}{2}\frac{1}{K-1}\frac{1}{|B_j|}\frac{1}{m_{ij}(n_in_j-m_{ij})};\\
    &\mathbb{P}_{(B_j,v,\vec{c}_v)}(W_1=(B_i,u,\vec{c}_u))=\frac{1}{2}\frac{1}{K-1}\frac{1}{|B_i|}\frac{1}{m_{ij}(n_in_j-m_{ij})},
\end{align*}
and detailed balance holds.
\end{itemize}

\subsection{Proof of Lemma \ref{lem:covermix2}}
\label{app:covermix2}
Based on the edgelighter walk in $V$, consider the chain $W_t=(B_t,L_t,C_t)$ where $L_t$ denotes the position of the edgelighter in $V$ at time $t$, $B_t$ denotes the community of $L_t$ and $C_t \in \{0,1\}^{\binom{n}{2}}$ the edge configuration at time $t$. 
Let $\tilde t^{\bullet}_{\text{cov}}$ be the first time that the walker has traversed every edge/non-edge within each community (i.e., for each $\{i,j\}\in\binom{V}{2}$ such that $i,j$ belong to the same community, the walker has made the $i \rightarrow j$ transition or/and the $j \rightarrow i$ transition). 

\noindent Let $\eta$ be an initial distribution for $W_0$ such that
\begin{itemize}
\item[i] The community $B_0$ is chosen uniformly randomly; 
\item[ii] Given $B_0=b$, $L_0$ is uniformly distributed over $b$;
\item[iii] For the initial edge configuration $C_0$, all edges within each community are absent, and all edges across communities are independently uniformly random (i.e., if $C_{0}^{i,j}$ is the portion of edges across blocks $i$ and $j$, then the $C_{0}^{i,j}$'s are independent and $\mathbb{P}(C_{0}^{i,j}=c_0^{(i,j})=1/\binom{n_in_j}{m_{ij}}$).
\end{itemize}
Let $c'$ be any fixed configuration such that all edges within each community are present in $c'$.
Given $\tilde t^{\bullet}_{\text{cov}}\leq t$, for $W_t=(b,x,c')$, all edges within each community need to have been flipped by time $t$ (this has probability $\prod_{j}q_{j;2}^{\binom{n_j}{2}}$); and given the community $b$, $x$ is uniformly distributed over $b$.
Such a $W_t$ may not put uniform mass on each community, but there exists a community $b$ and vertex $y\in b$ such that 
\begin{align*}
\mathbb{P}_\eta(W_t=(b,x,c'))
&=\mathbb{P}_\eta(W_t=(b,x,c')|\tilde t^{\bullet}_{\text{cov}}\leq t)\mathbb{P}(\tilde t^{\bullet}_{\text{cov}}\leq t)\\
&\leq \mathbb{P}_\eta(\tilde t^{\bullet}_{\text{cov}}\leq t) \frac{1}{|K|}\frac{1}{|b|}\prod_{j}q_{j;2}^{\binom{n_j}{2}}\frac{1}{\prod_{i,j}\binom{n_in_j}{m_{ij}}} .
\end{align*}
From this, we have that 
\begin{align*}
    1-\frac{\mathbb{P}_\eta(W_t=(b,x,c'))}{\pi^\bullet(b,x,c')}\geq 1-\mathbb{P}_\eta(\tilde t^{\bullet}_{\text{cov}}\leq t) .
\end{align*}
Let $s^\bullet(t,\eta)=\max_{(b,y,c)}1-\frac{\mathbb{P}_\eta(W_t=(b,y,c))}{\pi^\bullet(b,y,c)}$, and note that 
$\mathbb{P}(\tilde t^{\bullet}_{\text{cov}}> t) \leq s^\bullet(t,\eta)$.
Next, define
\begin{align*}
d^\bullet(t,\eta)&=\|\mathbb{P}^t_{\eta}(\cdot)-\pi^\bullet\|_{TV}\\ 
d^\bullet(t)&=\max_{(b,x,c)}\|\mathbb{P}^t_{(b,x,c)}(\cdot)-\pi^\bullet\|_{TV}\\ 
\bar d^\bullet(t,\eta)&=\max_{(b,y,c)}\|\mathbb{P}^t_{\eta}(\cdot)-\mathbb{P}^t_{(b,y,c)}(\cdot)\|_{TV}~~~\text{ and }\\
\bar d^\bullet(t)&=\max_{(b,x,c),(b',y,c')}\|\mathbb{P}^t_{(b,x,c)}(\cdot)-\mathbb{P}^t_{(b',y,c')}(\cdot)\|_{TV} .
\end{align*}
As in the proof of Lemma \ref{lem:covermix}, we have that, as $(Z_t)$ is reversible,
\begin{equation*}
d^\bullet(t,\eta)\leq \bar d^\bullet(t,\eta)\leq d^\bullet(t,\eta)+d^\bullet(t)
~~~\text{ and}~~~
s^\bullet(2t,\eta)\leq 1-(1-\bar d^\bullet(t,\eta))^2,
\end{equation*}
where the latter inequality controlling $s^\bullet(2t,\eta)$ follows as in the proof of Lemma 4.7 in \cite{aldous2002reversible}.
While $\bar d^\bullet(t,\eta)$ may not be submultaplicative, we do have that $\bar d^\bullet(t+s,\eta)\leq \bar d^\bullet(s,\eta)\bar d^\bullet(t)$.
This follows from the following argument, adapted from the proof of Lemma 4.11 in \cite{levin2017markov}:
Fix $(b,x,c)$ and let $(U_s,V_s)$ be the optimal coupling of $\mathbb{P}^s_{(b,x,c)}(\cdot)$ and $\mathbb{P}^s_{\eta}(\cdot)$.
We then have for any set $A$,
\begin{align*}
    \mathbb{P}^{s+t}_{(b,x,c)}(A)-\mathbb{P}^{s+t}_{\eta}(A)
    &=\mathbb{E}(\mathbb{P}^{t}(U_s,A)-\mathbb{P}^{t}(V_s,A))\\
    &\leq\mathbb{E}(\mathds{1}\{U_s\neq V_s)\bar d^\bullet(t))=\bar d^\bullet(t)\mathbb{P}(U_s\neq V_s)\leq \bar d^\bullet(s,\eta)\bar d^\bullet(t),
\end{align*}
as desired.

Let $t^\bullet(\cdot)$ be the usual mixing time of $W_t$.
As in the proof of Lemma \ref{lem:covermix}, we then have that, noting that the chain is reversible,
$$\bar d^\bullet(2t^\bullet_m,\eta)\leq
\bar d^\bullet(t^\bullet_m,\eta)\bar d^\bullet(t^\bullet_m)
\leq 2  (d^\bullet(t^\bullet_m,\eta)+d^\bullet(t^\bullet_m)) d^\bullet(t^\bullet_m)\leq 1/4.$$
and so,
$$\mathbb{P}_\eta(\tilde t^{\bullet}_{\text{cov}}> 4t^\bullet_m)\leq s^\bullet(4t^\bullet_m,\eta)\leq 1-(1-\bar d^\bullet(2t^\bullet_m,\eta))^2<1/2.$$
As in the proof of Lemma \ref{lem:covermix}, we then have that 
$$
\mathbb{E}_\eta(\tilde t^{\bullet}_{\text{cov}})\leq 8t^\bullet_m.
$$
Regardless of the starting location, by time $\tilde t^{\bullet}_{\text{cov}}$, $W_t$ needs to have traversed all $\Theta(n^{2})$ edges/non-edges within the largest community.
This provides that 
$$
\Omega(n^2)=\mathbb{E}_\eta(\tilde t^{\bullet}_{\text{cov}})\leq 8t^\bullet_m
$$
as desired.

\newpage
\subsection{Additional Experiments and Plots}
\subsubsection{Additional Standard edgelighter on ER model plots}
\label{app:ER_plt}

\begin{figure}[H]
    \centering
    \includegraphics[width=\textwidth]{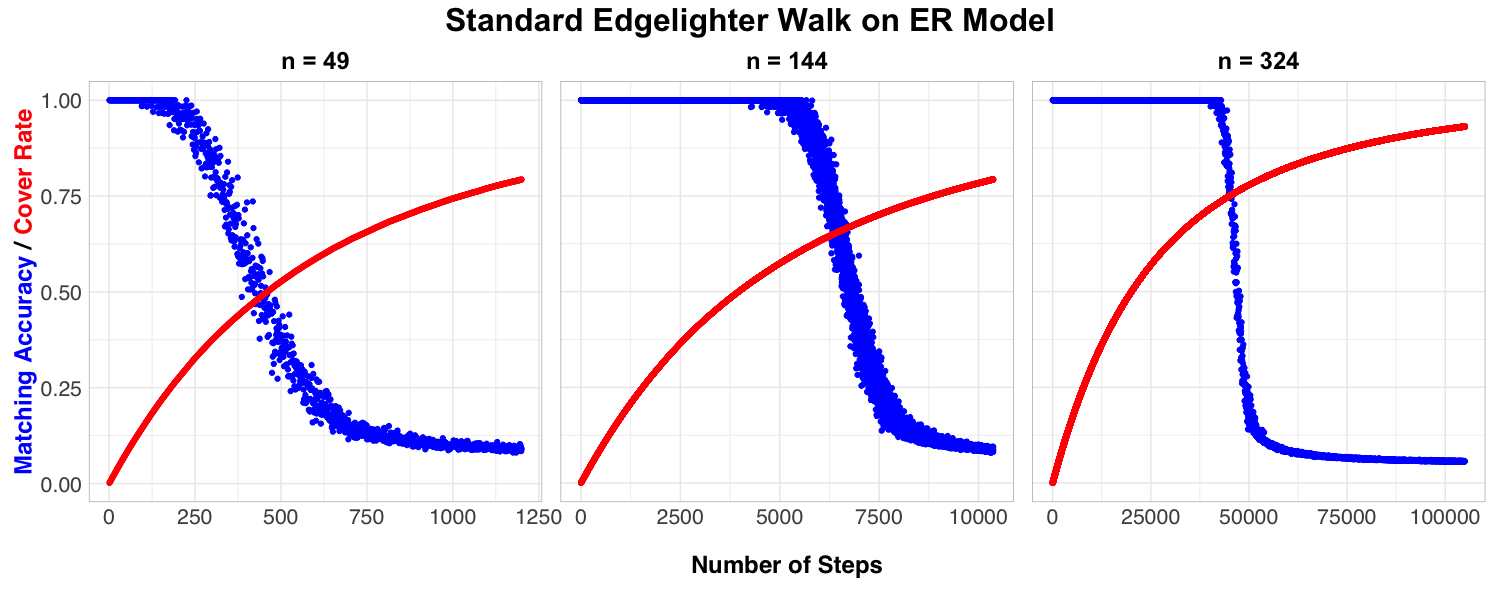}
    \caption{We plot matching correctness vs. number of steps for $n=49$ (left), $n=144$ (middle) and $n=324$ (right), in blue, where the setup of the model and the parameters chosen are the same as in Section \ref{exp:ER}. For these plots, we further impose the cover rate of the edges vs number of steps, in red.}
    \label{app:global_ratio}
\end{figure}

\begin{figure}[H]
    \centering
    \includegraphics[width=\textwidth]{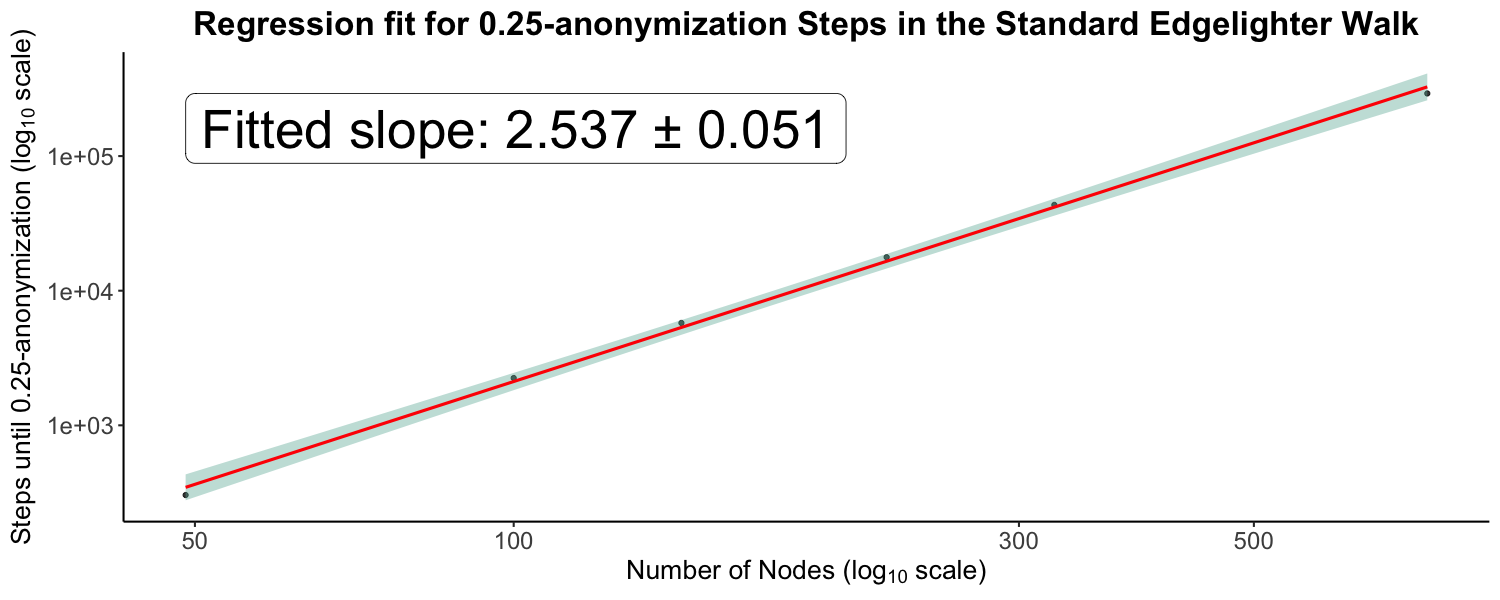}
    \includegraphics[width=\textwidth]{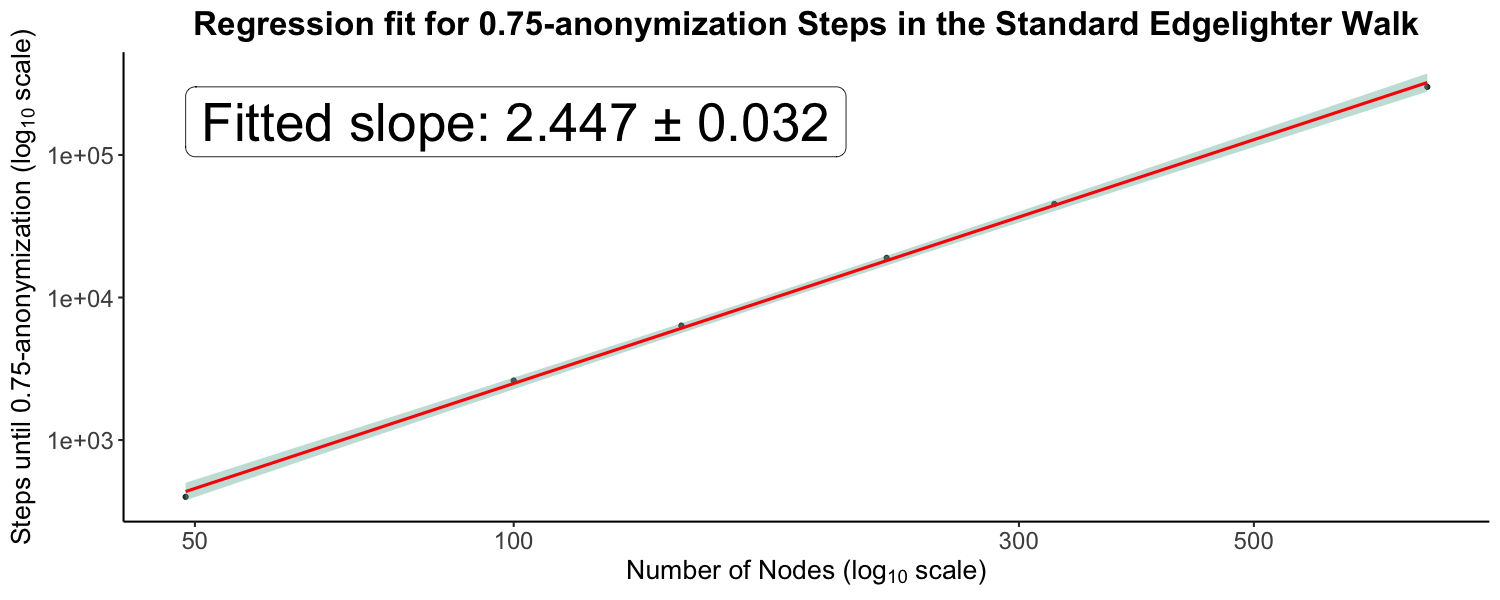}
    \caption{We plot a regression of number of steps needed until a {$0.25$}-anonymization (top panel) and number of steps needed until a {$0.75$}-anonymization (bottom panel) as a function of the number of nodes $n$ on a log-log scale. The setup of the model and the parameters chosen are the same as in Section \ref{exp:ER}}
\end{figure}

\subsubsection{Additional block edgelighter on SBM model plots}
\label{app:SBM_plt}

\begin{figure}[H]
    \centering
    \includegraphics[width=\textwidth]{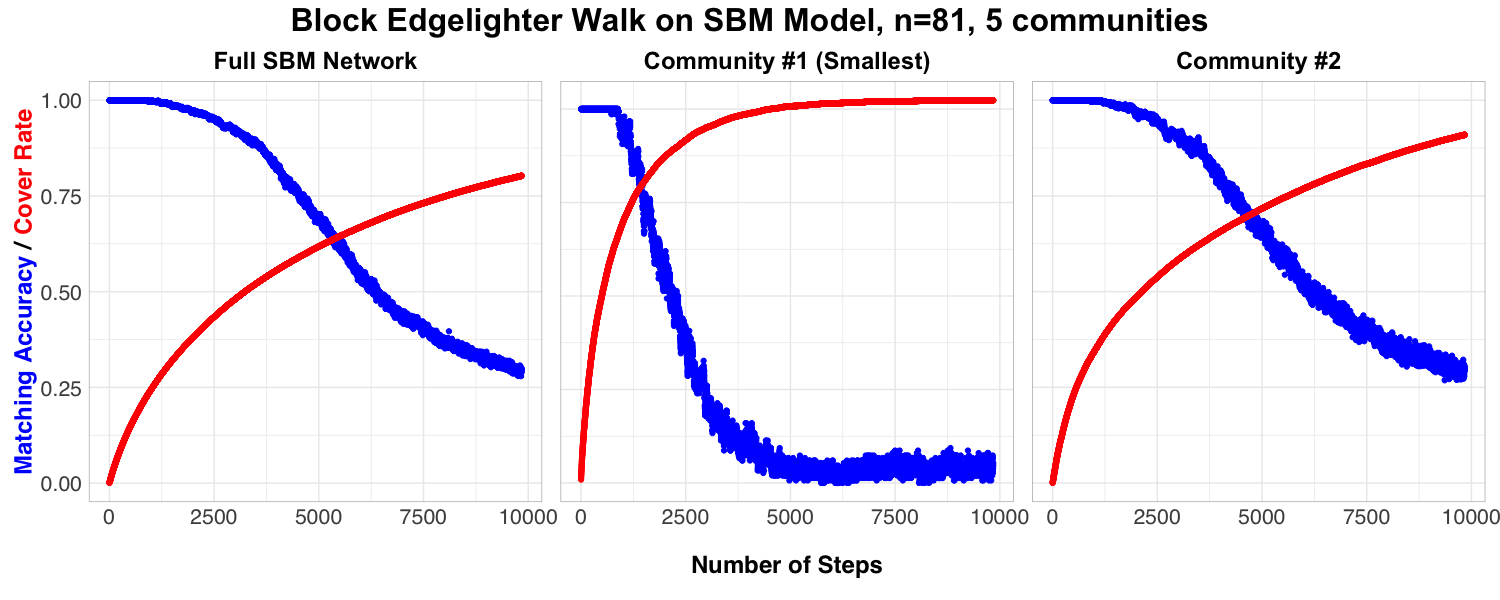}
    \includegraphics[width=\textwidth]{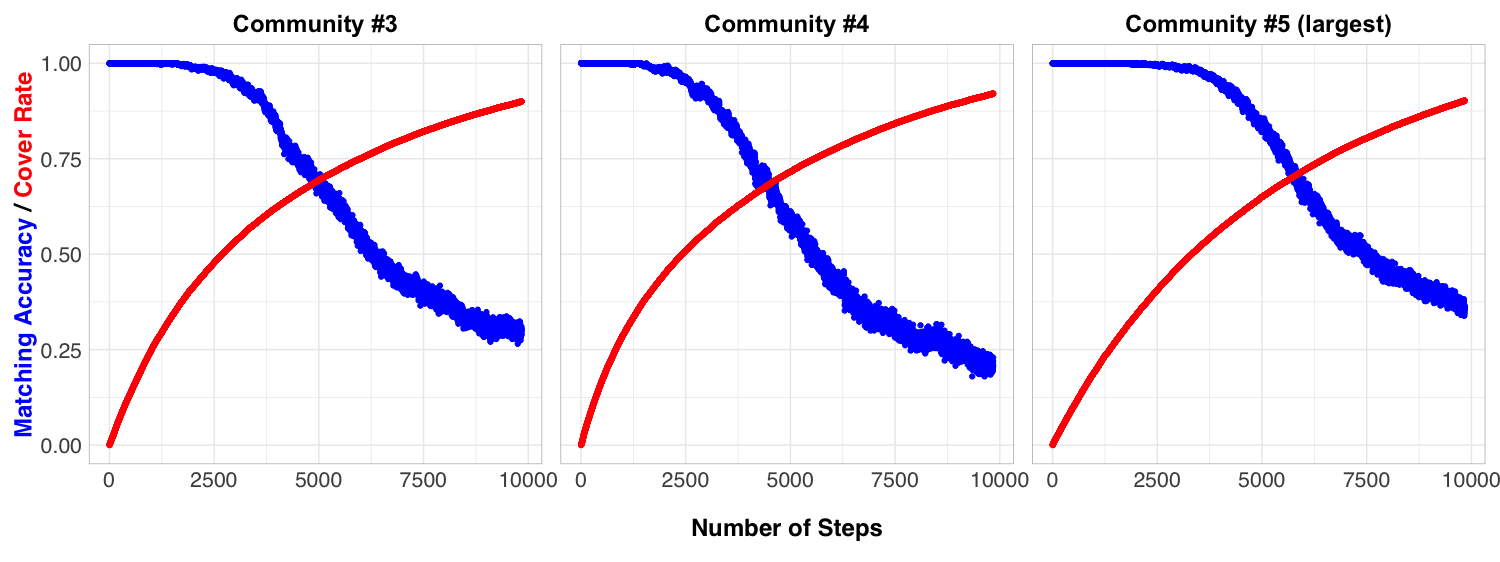}
    \caption{SBM edgelighter Walk plot with $n=81$, see Section \ref{exp:ER} for the description of the experiment setup. Number of communities $=5$}
    \label{fig:app_SBM_81}
\end{figure}

\begin{figure}[H]
    \centering
    \includegraphics[width=\textwidth]{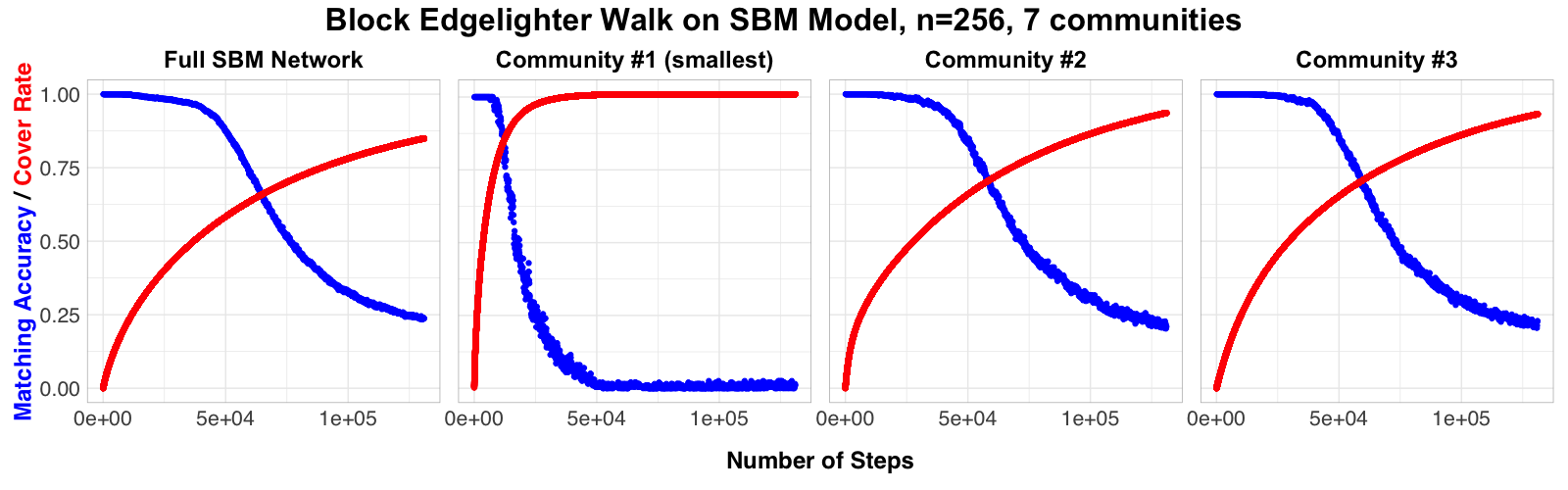}
    \includegraphics[width=\textwidth]{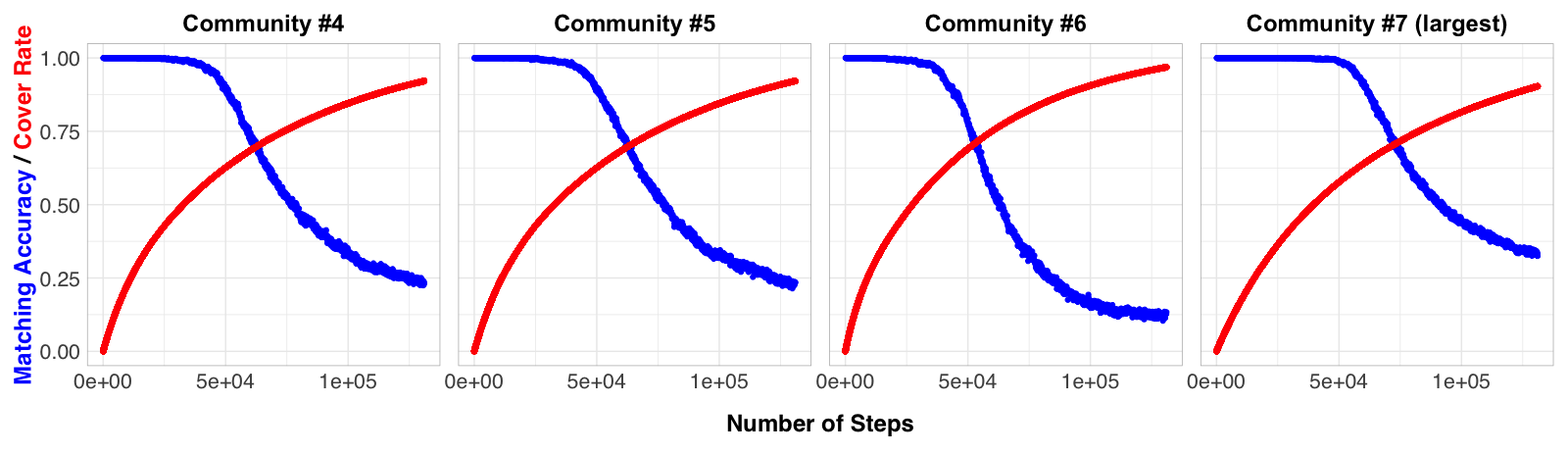}
    \caption{SBM edgelighter Walk plot with $n=256$, see Section \ref{exp:ER} for the description of the experiment setup. Number of communities $=7$}
    \label{fig:app_SBM_256}
\end{figure}

\begin{figure}[H]
    \centering
    \includegraphics[width=0.75\textwidth]{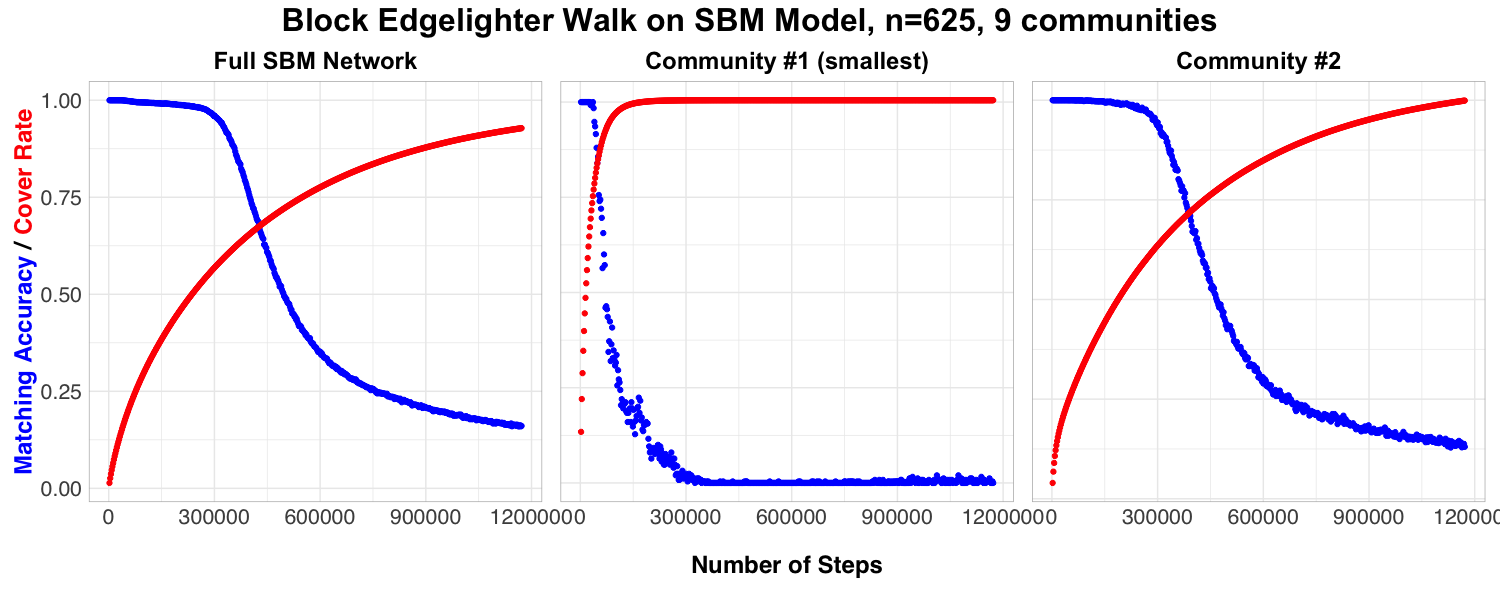}
    \includegraphics[width=0.75\textwidth]{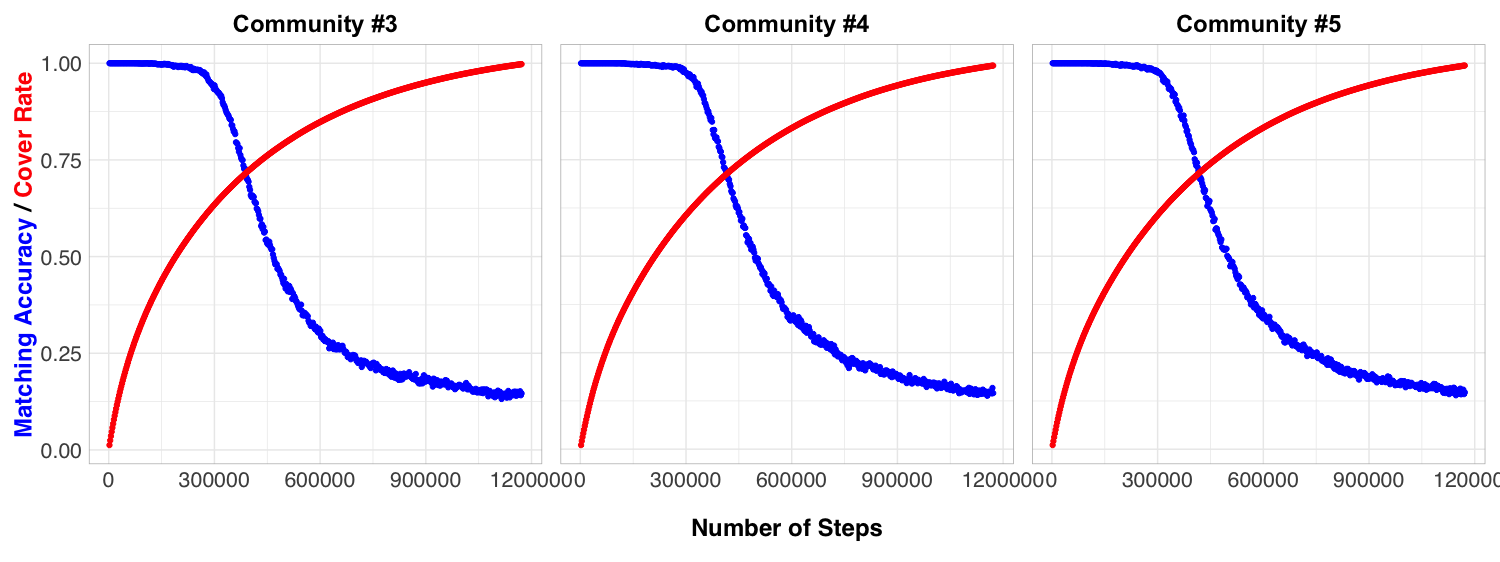}
    \includegraphics[width=\textwidth]{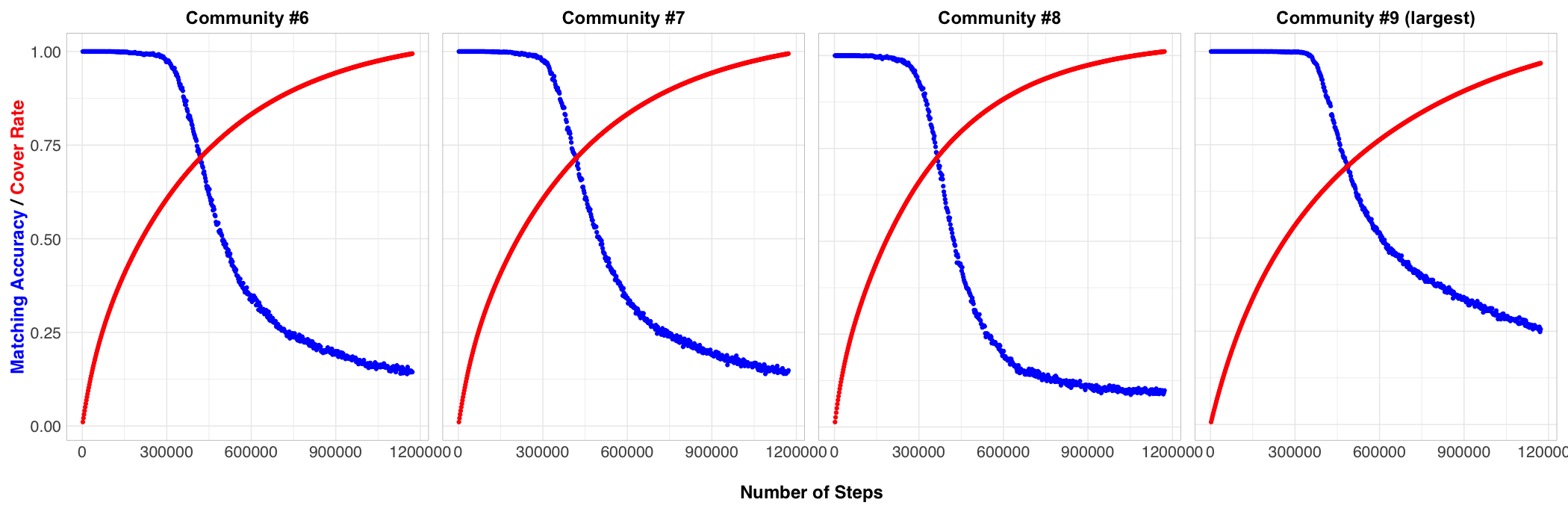}
    \caption{SBM edgelighter Walk plot with $n=625$, see Section \ref{exp:ER} for the description of the experiment setup. Number of communities $=9$}
    \label{fig:app_SBM_625}
\end{figure}

\begin{figure}[H]
    \centering
    \includegraphics[width=\textwidth]{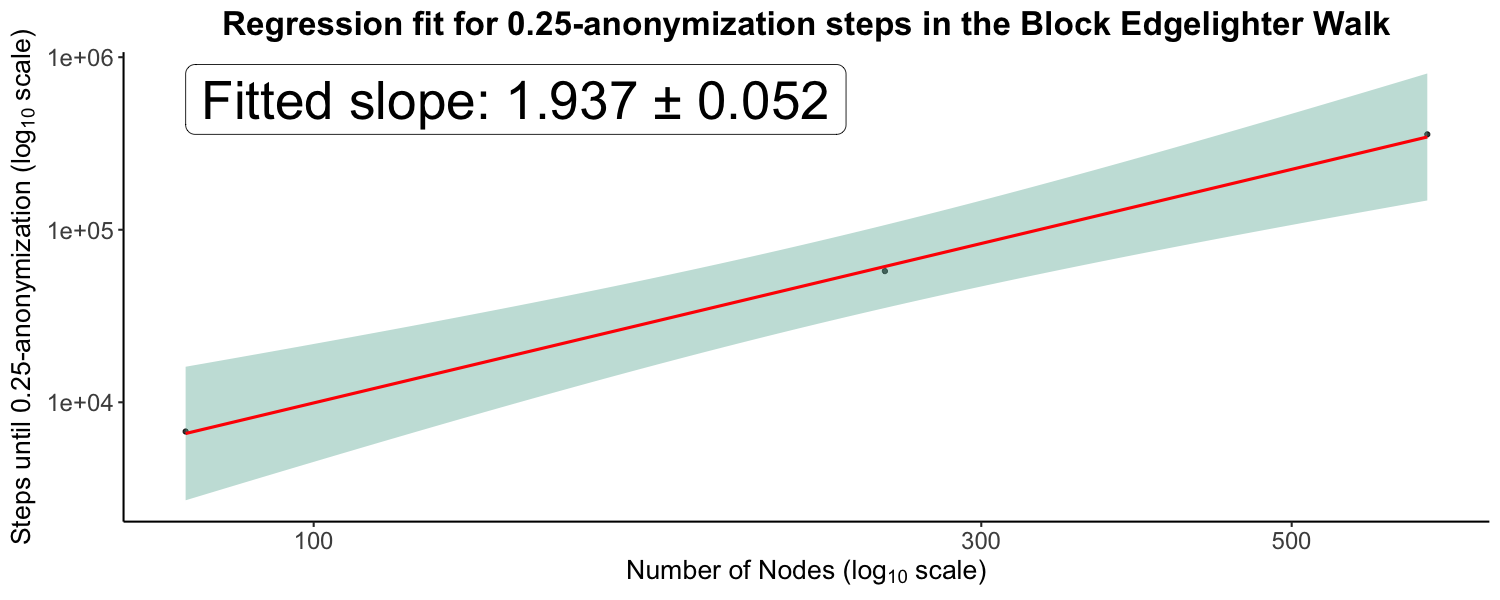}
    \includegraphics[width=\textwidth]{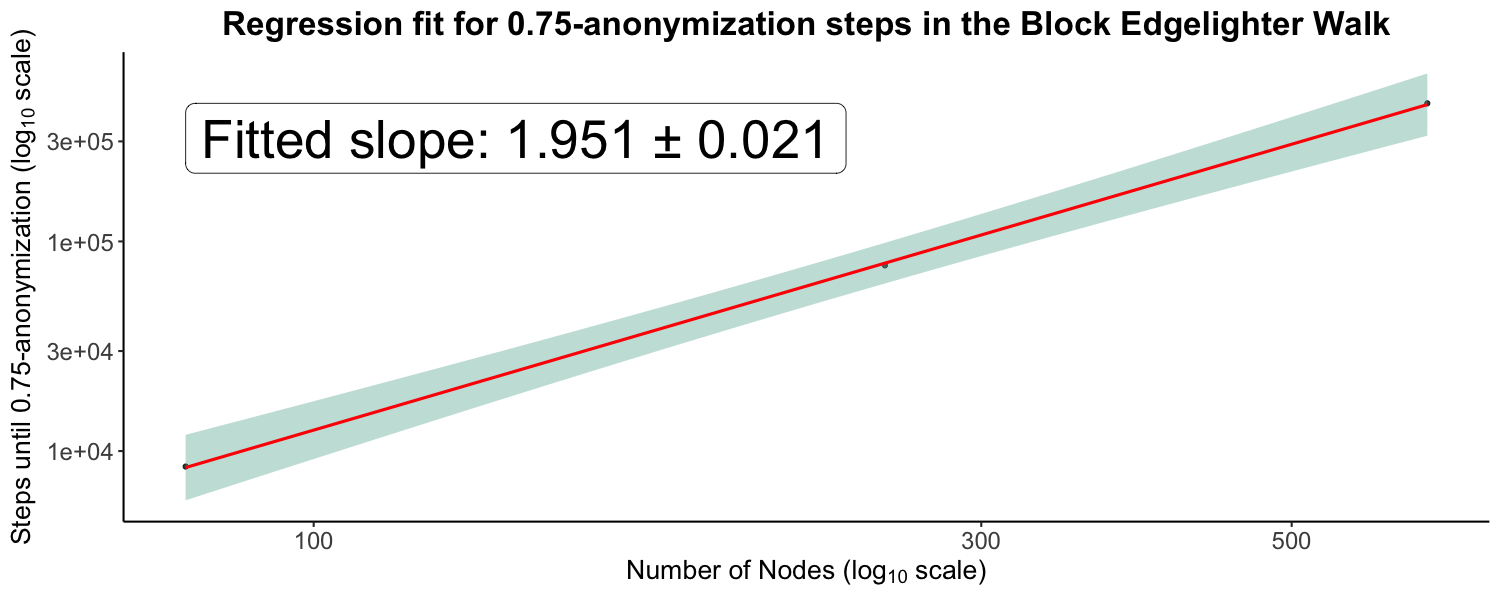}
    \caption{We plot a regression of number of steps needed until a {$0.25$}-anonymization (top panel) and number of steps needed until a {$0.75$}-anonymization (bottom panel) as a function of the number of nodes $n$ on a log-log scale. The setup of the model and the parameters chosen are the same as in Section \ref{exp:ER}}.
\end{figure}

\subsubsection{Additional block edgelighter on EU email network plots}
\label{app:EU_plt}

\begin{figure}[H]
    \centering
    \includegraphics[width=\textwidth]{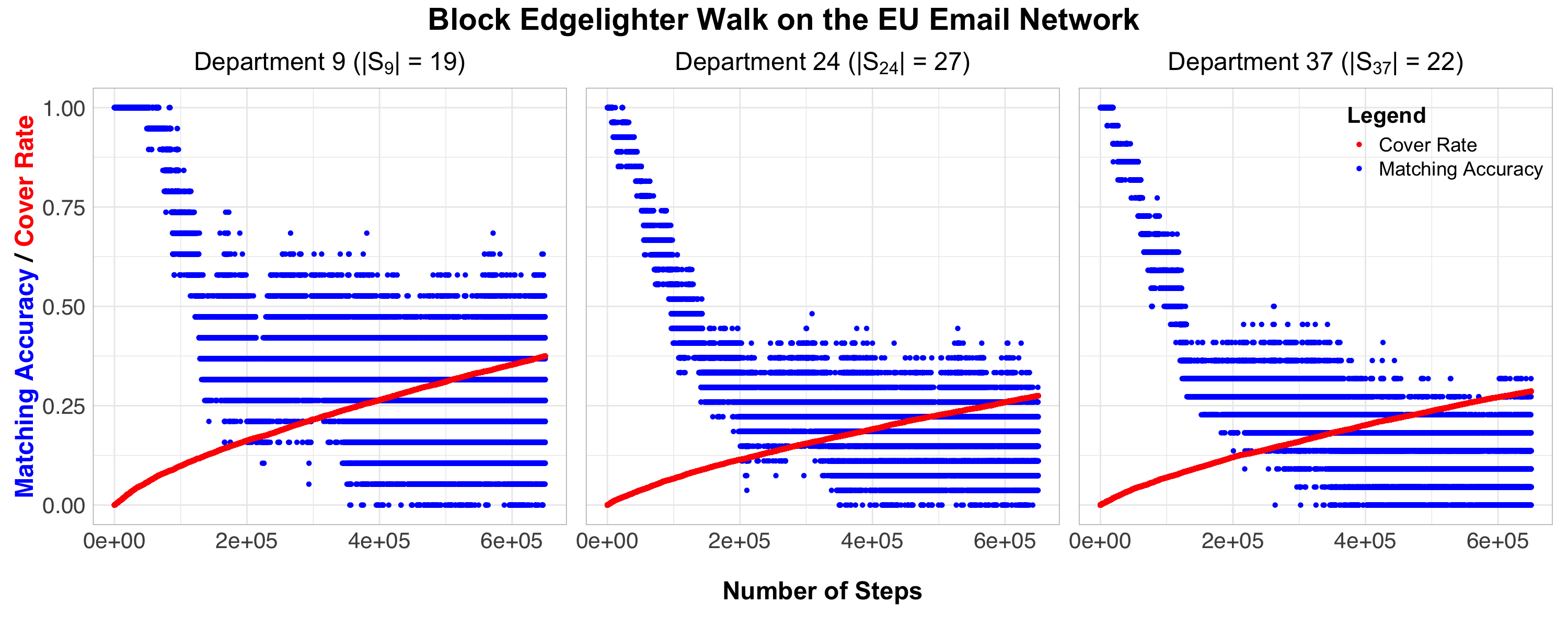}
    \includegraphics[width=\textwidth]{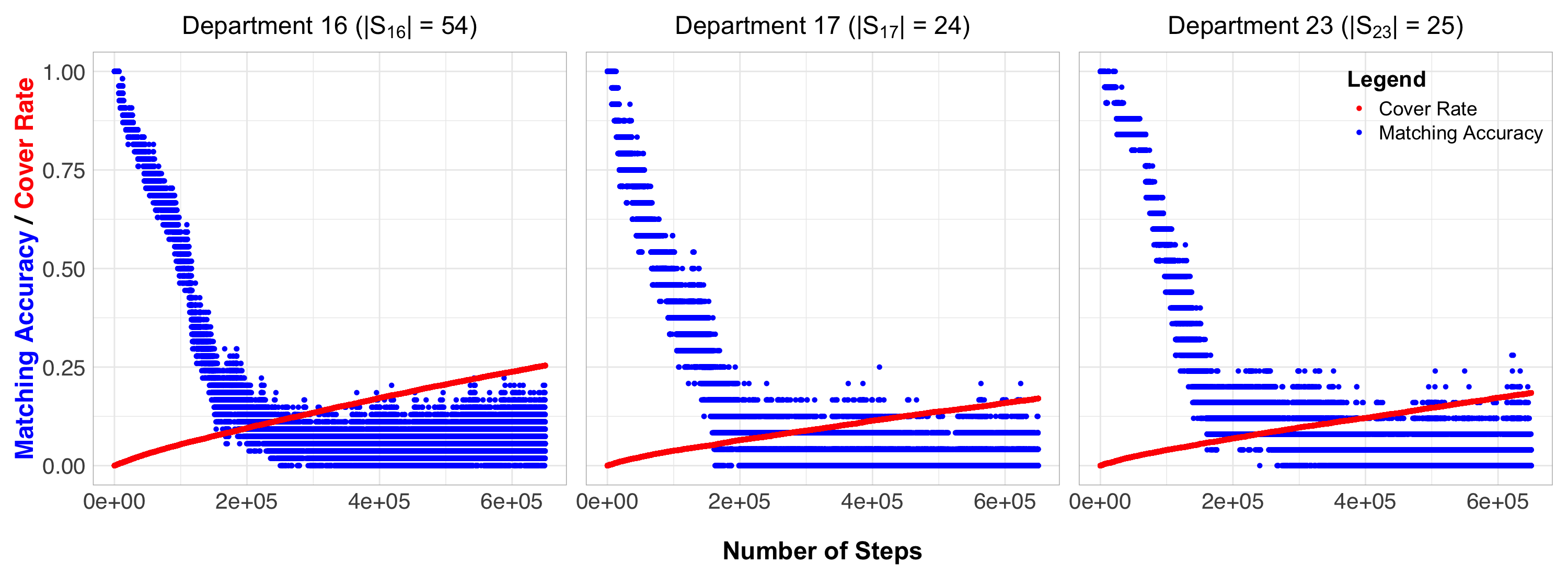}
    \caption{Block edgelighter walk plot for EU email network data, see Section \ref{sec:data} for the description of the experiment setup.}
\end{figure}
\end{document}